\documentclass[11pt, a4paper,leqno]{amsart}

\usepackage[centertags]{amsmath}
\usepackage{amsfonts}
\usepackage{amssymb}
\usepackage{mathtools}
\usepackage{mathrsfs}
\usepackage{amsthm}
\usepackage[english]{babel}
\usepackage[active]{srcltx}
\usepackage{enumerate}
\usepackage{esint}
\usepackage{xcolor}
\usepackage{graphicx}%
\usepackage[colorlinks=true, linkcolor=blue, citecolor=red, urlcolor=red]{hyperref}
\usepackage{bbm}
\usepackage{tikz}
\usepackage{float}
\usetikzlibrary{patterns}

\calclayout
\allowdisplaybreaks

\numberwithin{equation}{section}
\newtheorem{thm}{Theorem}[section]

\newtheorem{op}[thm]{Open Problem}
\newtheorem{defn}[thm]{Definition}

\def\N{{{\Bbb N}}}

\def\T{{{\Bbb T}}}

\def\R{\mathbb{R}}
\theoremstyle{remark}

\begin{document}

\title[On the DiPerna-Majda $2$D gap problem]{A sparse resolution of the DiPerna-Majda gap problem for $2$D Euler equations}

\author{Oscar Dom\'inguez}

\address{Oscar Dom\'inguez\\
	Departamento de M\'etodos Cuantitativos
	\\
	CUNEF Universidad
	\\
	28040 Madrid, Spain.}
\email{oscar.dominguez@cunef.edu}

\author{Daniel Spector}

\address{Daniel Spector\\ Department of Mathematics, National Taiwan Normal University, No. 88, Section 4,
Tingzhou Road, Wenshan District, Taipei City, Taiwan 116, R.O.C. \\
National Center for Theoretical Sciences\\No. 1 Sec. 4 Roosevelt Rd., National Taiwan
University\\Taipei, 106, Taiwan}
\email{spectda@protonmail.com}

\begin{abstract}
A central question which originates in the celebrated work in the 1980's of DiPerna and Majda asks what is the optimal decay $f > 0$ such that uniform rates $|\omega|(Q) \leq f(|Q|)$ of the vorticity maximal functions guarantee strong convergence without concentrations of approximate solutions to energy-conserving weak solutions of the $2$D Euler equations with vortex sheet initial data. A famous result of Majda (1993) shows $f(r) = [\log (1/r)]^{-1/2}$, $r<1/2$,
as the optimal decay for  \emph{distinguished} sign vortex sheets. In the general setting of \emph{mixed} sign vortex sheets,  DiPerna and Majda (1987)  established $f(r) = [\log (1/r)]^{-\alpha}$ with $\alpha > 1$ as a sufficient condition for the lack of concentrations, while the expected gap $\alpha \in (1/2, 1]$ remains as an open question. In this paper we resolve the DiPerna-Majda $2$D gap problem: In striking contrast to the well-known case of distinguished sign vortex sheets,  we identify $f(r) = [\log (1/r)]^{-1}$ as the optimal regularity for mixed sign vortex sheets that rules out concentrations.

For the proof, we propose a novel method to construct explicitly  solutions with mixed sign to the $2$D Euler equations in such a way that wild behaviour creates within the relevant geometry of \emph{sparse} cubes (i.e., these cubes are not necessarily pairwise disjoint, but their possible overlappings can be controlled in a sharp fashion). Such a strategy is inspired by the recent work of the first author and Milman \cite{DM} where  strong connections between energy conservation and sparseness are established.

\end{abstract}

\maketitle

\section{Introduction}

In two dimensions, the Euler equations for an inviscid incompressible fluid flow are given by
\begin{equation}
\left\{
\begin{array}
[c]{c}%
\partial_t u+u\cdot\nabla u=-\nabla p,\\
\text{div }u=0,\\
u(0, \cdot) = u_0,
\end{array}
\right.  \label{Euler}%
\end{equation}
where $u:\mathbb{R}^2 \times [0, T] \to \mathbb{R}^2$ denotes the  \emph{fluid velocity}, $p:\mathbb{R}^2 \times [0, T] \to \mathbb{R}$ the (scalar) \emph{pressure}, and $u_0:\mathbb{R}^2 \to \mathbb{R}^2$ the \emph{initial fluid velocity}.  Further denote by $\omega:= \operatorname*{curl}u:\mathbb{R}^2 \times [0, T] \to \mathbb{R}$ the \emph{vorticity} and $\omega_0:= \operatorname*{curl}u_0:\mathbb{R}^2  \to \mathbb{R}$ the \emph{initial vorticity}, two objects of central interest in the study of these equations.

Motivated by approximations arising in numerical implementations and formal asymptotic analysis of Euler equations  (e.g. vortex blob approximation, vanishing viscosity limits of the Navier-Stokes equations, and smoothing of initial data for exact solutions of the Euler equations), in a series of influential papers in the late 80's, DiPerna and Majda \cite{DiPernaMajda, DiPernaMajda1*, DiPernaMajda2} developed  a rigorous framework of approximate solutions to \eqref{Euler} for \emph{vortex sheet} initial data
\begin{equation}\label{Vortex}
	\omega_0 \in \mathcal{M} \cap H^{-1}_{loc},
\end{equation}
where $\mathcal{M}$ is the space of Radon measures.  We recall that $\{u^{\varepsilon}\}, \, \varepsilon > 0,$ is said to be an \emph{approximate solution family} of \eqref{Euler} if the following properties are satisfied:
\begin{enumerate}
	\item[(P1)] $\{u^{\varepsilon}\}$ is uniformly bounded in $L^{\infty}([0,T];L^{2}_{loc})\cap\text{Lip}%
((0,T);H^{-L}_{loc})$ for some\footnote{The
uniform bound in $\text{Lip}((0,T); H^{-L}_{loc})$ is a technical assumption in order to guarantee that
initial vector fields $u^{\varepsilon}(0, \cdot)$ are well-defined. In
practice, this follows easily from the uniform energy bound $L^{\infty
}([0,T];L^{2}_{loc})$.} $L>1$.
\item[(P2)]  Weak consistency with \eqref{Euler} in the sense that
\begin{equation}\label{P2a}
\int_{0}^{T}\int_{\mathbb{R}^{2}}\varphi_{t}\cdot u^{\varepsilon}%
+(D\varphi\,u^{\varepsilon})\cdot u^{\varepsilon}\,dx\,dt+\int_{\mathbb{R}%
^{2}}\varphi(x,0)\cdot u^{\varepsilon}(x,0)\,dx\rightarrow0\qquad
\text{as}\quad\varepsilon\rightarrow0 
\end{equation}
for every test field $\varphi\in C^{\infty}$ with $\operatorname{div}\varphi=0$. Here $D\varphi$ is
the Jacobian matrix of $\varphi$. 
\item[(P3)] $\operatorname{div}u^{\varepsilon}=0$ (in the distributional sense).
\end{enumerate}

Observe that the energy bound (P1) guarantees that (possibly passing to a subfamily) $u^\varepsilon \rightharpoonup u$ in $L^\infty([0, T]; L^2_{loc})$.  The question of interest is then to obtain necessary and sufficient conditions for which the weak limit $u$ is a weak solution of \eqref{Euler}.  If there is no concentration, i.e. one is able to upgrade the weak convergence $u^\varepsilon \rightharpoonup u$ to a strong convergence $u^\varepsilon \to u$ in $L^\infty([0, T]; L^2_{loc})$, then one has the convergence of  quadratic terms $(D\varphi\,u^{\varepsilon})\cdot u^{\varepsilon}$ in \eqref{P2a} and hence $u$ will be a weak solution.  This is the case in a well-known example, when $\omega_0 \in L^p_c, \, p > 1$, while one of the main results of  the DiPerna--Majda theory relies on $\log$-Morrey norms\footnote{More precisely, the authors refers to \eqref{Morrey} as \emph{uniform decay of vorticity maximal function.}}:
\begin{equation}\label{Morrey}
	\|\omega\|_{M^{\log, \alpha}} = \sup_Q \, (1 - (\log |Q|)_-)^{\alpha} |\omega|(Q), \qquad \alpha > 0, 
\end{equation}
where\footnote{A simple consequence of H\"older's inequality gives $L^p \subset M^{\log, \alpha}$ for any $p > 1$ and $\alpha > 0$.} the supremum runs over all cubes $Q$ with sides parallel to the axes of coordinates and $|\omega|(Q)$ is the total variation of $\omega \in \mathcal{M}$ in $Q$. The finiteness of \eqref{Morrey} admits an interpretation in terms of circulation decays at uniform rates created in small scales. The following result was first established  in \cite[Theorem 3.1]{DiPernaMajda}, while several alternative (and simplified) proofs were later obtained in \cite{LNT, CS, J, Lan2, DM}.

\begin{thm}\label{ThmDMConser}
Assume that $\{u^\varepsilon\}$ is an approximation family of the $2$D Euler equations \eqref{Euler} such that the corresponding family of vorticities
\begin{equation}\label{HDM1}
 \{\omega^\varepsilon\} \quad  \text{is uniformly bounded in} \quad L^\infty ([0, T]; M^{\log, \alpha}) \quad \text{for some} \quad \alpha > 1.
\end{equation}
Then $\{\omega^\varepsilon\}$ defines a bounded family of vortex sheets  \eqref{Vortex} and\footnote{Throughout this paper, \eqref{HDM2} should be interpreted modulo passing to a subfamily, if necessary.} 
\begin{equation}\label{HDM2}
\lim_{\varepsilon \to 0} u^\varepsilon = u \quad \text{in} \quad L^\infty ([0, T]; L^2_{loc}),
\end{equation}
where $u$ is a weak solution to \eqref{Euler}. 
\end{thm}

Moreover it is known that, still under the assumption \eqref{HDM1} for relevant choices of approximate solutions, the conclusion of  Theorem \ref{ThmDMConser} can be strengthened to say that the weak solution $u$ conserves the kinetic energy; see Theorem \ref{ThmConser} below.  However, another remarkable scenario can arise, when the weak convergence $u^\varepsilon \rightharpoonup u$ fails to be strong and yet still $u$ is a weak solution of \eqref{Euler}, so-called \emph{concentration-cancellation} phenomenon.  In sharp contrast with solutions obtained from  Theorem \ref{ThmDMConser}, solutions arising from concentration-cancellation exhibit in general a quite ``wild" behaviour and, in particular, they do not preserve energy.   Hence there exists  a strong connection between lack of concentrations vs. concentration-cancellation and the important dichotomy energy conservation vs. anomalous dissipation. See Section \ref{SectionECAD} below for further details. 

According to \cite[Remark 3.2, p.~337]{DiPernaMajda}, Theorem \ref{ThmDMConser} is \emph{nearly} optimal in the sense that explicit examples of vorticities  were constructed in \cite{DiPernaMajda1*} exhibiting concentration-cancellation  in the limit process when the regularity assumption $M^{\log, \alpha}, \, \alpha > 1,$ in \eqref{HDM1} is relaxed\footnote{Note that $M^{\log, \alpha_0} \subset M^{\log, \alpha_1}$ if $\alpha_0 > \alpha_1$.} to  $M^{\log, 1/2}$.  Interestingly, this optimality assertion has its roots on the celebrated result of Delort \cite{Delort} (more precisely, its alternative proof due to Majda \cite{M93}) showing global existence of weak solutions to \eqref{Euler} for vortex sheets initial data with distinguished sign (i.e., $\mathcal{M}$ is replaced by $\mathcal{M}_+$ in \eqref{Vortex}, where $\mathcal{M}_+$ is the space of nonnegative Radon measures). Moreover,  the fact \cite[Proposition, p. 928]{M93} that $\{\omega^\varepsilon\}$ is uniformly bounded in $M^{\log, 1/2}$ for $\{\omega^\varepsilon\} \subset \mathcal{M}_+ \cap H^{-1}_{loc}$ suggests that $M^{\log, 1/2}$ is the borderline regularity space for concentration-cancellation, at least when  attention is restricted to vortex sheets with distinguished sign. However, going from distinguished to mixed sign vortex sheets is a major open problem, which finds many physical motivations and higher complexity generated by irregular flows intertwining positive and negative regions as illustrated by  the numerical performances from \cite{K}.  To close  the gap between regularity in  $M^{\log, \alpha}, \alpha > 1,$ guaranteeing strong convergence to weak solutions and  $M^{\log, 1/2}$-regularity  where concentrations arise remains as an outstanding open question for general vorticities. This can be informally formulated as:

\begin{op}[DiPerna-Majda gap problem, 1987]\label{OPDM}
	What is the optimal value of $\alpha \in [\frac{1}{2}, 1]$ for which concentration-cancellation phenomenon occurs in $M^{\log, \alpha}$? Or in other words, is it possible to extend the lack of concentrations given in Theorem \ref{ThmDMConser} to $M^{\log, \alpha}$ for some  $\alpha \in (\frac{1}{2}, 1]$?
\end{op}

Considerable amounts of effort have been undertaken to solve this question, in particular, we refer the reader to Section \ref{History} below to find an overview on the state of the art.

%%%%%%%%%%%%%%%%%%%%%

\subsection{The main results} \hfill\\

In this paper, we resolve  Open Problem \ref{OPDM}. Specifically, as a first step we  show that the vortex sheet problem is ill-posed in $M^{\log, 1}$ in the following sense:

\begin{thm}\label{ThmMain}
	There exists a positive-signed Radon probability measure $\omega$ (related to a certain Cantor set $E \subset [0, 1]^2$)  such that $\omega \in M^{\log, 1}$ but $\omega \not \in H^{-1}_{loc}$.
\end{thm}

This result tells us that distinguished sign vorticities in $M^{\log, 1}$ may generate irregular flows that do not define vortex sheets. This is in striking contrast with $M^{\log, \alpha}, \, \alpha > 1$, which is always formed by vortex sheets, see Theorem \ref{ThmDMConser}.

As a second step we show that  $M^{\log, 1}$ is the optimal regularity space for vortex sheets concentrations, in the sense that wild behaviour of vortex sheets may arise under this regularity assumption, but no concentration occurs in $M^{\log, \alpha}$ if $\alpha > 1$ (see Theorem \ref{ThmDMConser}). This gives an answer to the gap problem posed in  Open Problem \ref{OPDM} by substantially improving the best earlier known optimality assertions from $M^{\log, 1/2}$ (Delort's result) to the smaller space $M^{\log, 1}$. Moreover, we are able to characterize  concentration sets in terms of the measure $\omega$ constructed in Theorem \ref{ThmMain}. The precise statement reads as follows.

\begin{thm}[Concentration-cancellation in $M^{\log, 1}$]\label{ThmMain2}
	There exists a family $\{u^\varepsilon\}$ of exact steady solutions to \eqref{Euler} satisfying the following properties:
	\renewcommand\labelitemi{\tiny$\bullet$}
	\begin{itemize}
		\item $\{u^\varepsilon\}$ is bounded in $L^2_{loc}$ and $\{\omega^\varepsilon\}$ is a bounded set in   $M^{\log, 1}$ formed by mixed-sign vortex sheets;
		\item $\{u^\varepsilon\}$ converges weakly to the trivial solution of \eqref{Euler} but not strongly in $L^2$;
		\item the reduced defect measure $\theta$ for $\{u^\varepsilon\}$ (cf. \eqref{EnergyConcentration}) is comparable to the measure $\omega$ given in Theorem \ref{ThmMain}: For an open set $A \subset \R^2$, 
		$$
			A  \text{ is a concentration set for } \{u^\varepsilon\} \iff \omega(A) >0.
		$$  
		In particular, the Cantor set $E$ from Theorem  \ref{ThmMain} is a concentration set with zero Hausdorff dimension. 
	\end{itemize}
\end{thm}

Recall that $\theta :  \mathcal{B}(\R^2 \times \R^+) \to [0,\infty]$ is the \emph{reduced defect measure}\footnote{Note that $\theta$ is not a measure, and fails to be even countably subadditive.} defined by
\begin{equation}\label{EnergyConcentration}
	\theta(E)  := \limsup_{\varepsilon \to 0} \int_E |u^\varepsilon(x, t) - u(x, t)|^2 \, dx \, dt, 
\end{equation}
which precisely characterizes regions where weak convergence of $\{u^\varepsilon\}$ fails to be strong.  A detailed study of $\theta$ as a tool to measure weak convergence was developed  in \cite{DiPernaMajda2}, as it is useful to distinguish between lack of concentrations (and then existence of weak solutions to \eqref{Euler} that are obtained as strong limit of approximate solutions) and concentration-cancellation.

We point out that Theorem \ref{ThmMain2} does not contradict the above mentioned fact that $M^{\log, 1/2}$ is the optimal regularity space for concentration-cancellation when restricted to distinguished sign vortex sheets. Indeed, the corresponding family $\{\omega^\varepsilon\}$  in Theorem \ref{ThmMain2} is formed by mixed-sign vortex sheets.

\subsection{Earlier contributions on the gap problem}\label{History} \hfill\\

Next we review some previous results  related to Open Problem \ref{OPDM}, focusing only on the most relevant ones for the purposes of this paper, as well as its strong connections with the longstanding question of energy conservation and anomalous dissipation.  

\subsubsection{\emph{\textbf{Connections with energy conservation and anomalous dissipation}}}\label{SectionECAD}
It is plain to see that smooth solutions $u$ to \eqref{Euler} conserve the kinetic energy in the sense that\footnote{To avoid further technicalities, we switch temporarily from  $u \in L^2_{loc}(\R^2)$ to  $u \in L^2(\T^2)$, where $\T^2$ is the $2$-dimensional torus.}
\begin{equation}\label{Conservative}
	\|u(t)\|_{L^2(\T^2)} = \|u_0\|_{L^2(\T^2)} \qquad \forall t \in [0, T]. 
\end{equation}
By contrast, the validity of \eqref{Conservative} for \emph{weak} solutions is not guaranteed, leading to the famous anomalous dissipation phenomenon,  i.e., irregular inviscid flows that solve \eqref{Euler} but do not preserve the energy.  Anomalous dissipation is a cornerstone of turbulence theory and there is a huge literature on this topic, mainly related to the Onsager's conjecture 
 on the critical H\"older regularity in $3$D for energy conservation. We do not aim to provide here with a detailed discussion on the developments of Onsager's conjecture and its historical background and physical motivations, but  we only refer  to the surveys \cite{DlS2} and \cite{BV}.   However, we should at least mention that Onsager's conjecture was recently resolved by Isett \cite{I} (see also \cite{BDSV}): Anomalous dissipation may occur for weak solutions of \eqref{Euler} in the H\"older class $L^\infty([0, T]; C^{\alpha})$ with regularity $\alpha < 1/3$. We mention that a key role in the proof of \cite{I} is played by the celebrated convex integration method of De Lellis and Sz\'ekelyhidi  \cite{DlS}. Very recently, a proof of the  Onsager conjecture in $2$D has been obtained by Giri and Radu \cite{GR}.  On the other hand, conservation of energy \eqref{Conservative} is known to hold for weak solutions with regularity $C^\alpha, \, \alpha > 1/3$ (not only in $3$D, but also in $2$D), see  \cite{CET} and the extension to the borderline regularity $\alpha=1/3$ given in \cite{CCFS}. 
 
 It is well known that there exist  strong relationships between the phenomena of energy conservation/anomalous dissipation and strong convergence/concentration-cancellation of approximate solutions. Indeed, a classical result (cf. \cite{DiPernaMajda})  asserts that strong convergence of approximate solutions obtained via smoothing of the initial data yields solutions that preserve the energy \eqref{Conservative}. This result has been recently extended by Cheskidov, Lopes Filho, Nussenzveig Lopes, and Shvydkoy \cite{CLNS} to cover the important class of solutions obtained in the zero-viscosity limit, which are called    \emph{physically realizable weak solutions}. Namely, if $u$ is a physically realizable weak solution  with corresponding initial data $\omega_0 \in L^p, \, p > 1$, then \eqref{Conservative} holds. See also \cite{CCS}.  In particular, this statement highlights that Onsager's conjecture may not provide with definite answers on energy conservation since the optimal form of the rigid part of Onsager's conjecture \cite{CCFS} only applies to  $\omega_0 \in L^p$ with $p \geq 3/2$. In this vein, we also refer to the recent papers \cite{Lan} and \cite{DM}, where the role of $L^p$ in  the result of Cheskidov-Lopes Filho-Nussenzveig Lopes-Shvydkoy  is replaced by any $X$ rearrangement invariant and compactly embedded into $H^{-1}$, among other results. In particular, the following  statement that characterizes energy preservation may be found in \cite[Theorem 2.11]{Lan} (see also \cite[Theorem 2.8]{JLLN} for an extension to forced Euler equations.)
 
 \begin{thm}\label{ThmConser}
  Assume that $u \in L^\infty([0, T]; L^2)$ is a physically realizable weak solution to \eqref{Euler}. Then the following are equivalent:
 \begin{itemize}
 	\item $u^\varepsilon \to u$ in $C([0, T]; L^2)$,
 	\item $u^\varepsilon \to u$ in $L^p([0, T]; L^2), 1 \leq p < \infty$,
	\item $u$ conserves the energy \eqref{Conservative}.
 \end{itemize}
 \end{thm}
 
 In other words, the concentration-cancellation phenomenon for viscous approximations leads to anomalous dissipation.

\subsubsection{\emph{\textbf{The $H^{-1}$-stability method of Lopes Filho, Nussenzveig Lopes, and Tadmor.}}}\label{Sstab} A characterization of the lack of concentrations for approximate solutions $\{u^\varepsilon\}$ to Euler equations is provided by the powerful \emph{$H^{-1}$-stability method} \cite{LNT} of Lopes Filho, Nussenzveig Lopes, and Tadmor (with \cite[Section 4.2]{Lio96} as a forerunner). In particular, Theorem 1.1 from \cite{LNT} reads as follows.

\begin{thm}\label{THMLNT}
Let $\{u^\varepsilon\}$ be an approximation family of the $2$D Euler equations\footnote{In this paper, we are only interested in $2$D Euler equations, however we would like to mention that the criterion for lack of concentrations provided by Theorem \ref{THMLNT} also holds for higher dimensions.} \eqref{Euler} with related set of vorticities $\{\omega^\varepsilon\}$. 
 If $\{\omega^\varepsilon\}$ is uniformly bounded in a certain Banach function space  $X \subset \mathcal{M}$ (i.e., $\sup_{\varepsilon > 0} \sup_{t \in [0, T]} \|\omega^\varepsilon(t)\|_X < \infty$) such that
\begin{equation}\label{Compact}
	X \overset{compactly}{\hookrightarrow} H^{-1}_{loc},
\end{equation}
then 
\begin{equation*}
\lim_{\varepsilon \to 0} u^\varepsilon = u \qquad \text{in} \qquad L^\infty ([0, T]; L^2_{loc}),
\end{equation*}
where $u$ is a weak solution to \eqref{Euler}.  
\end{thm}

As a byproduct of Theorems \ref{ThmConser} and \ref{THMLNT}, physically realizable weak solutions with corresponding initial data $\omega_0 \in X$ for any rearrangement invariant space $X$ with \eqref{Compact}  preserve the energy  \eqref{Conservative}. Theorem \ref{THMLNT} provides a unified approach to identify critical regularity for lack of concentrations in wide classes of spaces $X$ including Lebesgue, Lorentz, Orlicz, and Morrey spaces. In particular,  Theorem \ref{ThmDMConser} can be recovered\footnote{In fact, the original proof of Theorem \ref{ThmDMConser} given in \cite{DiPernaMajda} requires an additional weak decay of vorticities at infinity, which can be overcome in the method of  \cite{LNT}.} from Theorem \ref{THMLNT}  since
\begin{equation}\label{CMorrey}
	M^{\log, \alpha}  \overset{compactly}{\hookrightarrow} H^{-1}_{loc} \qquad \text{if} \qquad \alpha > 1; 
\end{equation}
see \cite[Theorem 4.4]{LNT}. 
According to \cite[Theorem 4.2]{LNT}, the compactness assertion \eqref{CMorrey} is attributed to DeVore and Tao, independently.  

Note that the extension of \eqref{CMorrey} to $\alpha \in (1/2, 1]$ would be a sufficient condition to rule out concentration-cancellation in $M^{\log, \alpha}$ for a family of approximate solutions to the Euler equations, and in particular would resolve Open Problem \ref{OPDM}.  However, as mentioned explicitly by the authors in \cite[p. 400]{LNT}, the validity of such a compactness assertion was an intriguing open question.  In this paper, we resolve this question by showing that one does not even have the embedding  \eqref{CMorrey} for $\alpha \leq 1$, see \eqref{Dfv} and \eqref{Dfv2} below.

\subsubsection{\emph{\textbf{Tadmor's approach to Open Problem \ref{OPDM}: the role of packings}}}\label{SectionT}
Recall that a countable family of cubes $(Q_i)_{i \in I}$ is said to be a \emph{packing} if $Q_i \cap Q_j = \emptyset$ for $i \neq j$. Let us denote by $\Pi$ the set of all packings. Note that \eqref{Morrey} can be rewritten in terms of packings as
\begin{equation}\label{MorreyEquiv}
	\|\omega\|_{M^{\log, \alpha}} = \sup_{\Pi} \sup_{i \in I}  \, (1 - (\log |Q_i|)_-)^{\alpha} |\omega|(Q_i). 
\end{equation}
An interesting  attempt to solve Open Problem \ref{OPDM} was proposed by  Tadmor \cite{Tadmor}, who suggested to work with the new\footnote{We mention that $T^{\log, \alpha}$ is an special case of the functional classes $T^{p, q} \log^\alpha$, namely, $T^{\log, \alpha} = T^{1, 2} \log^\alpha$. Another distinguished example is  $T^{p, p} \log^0 = L^p$ (Riesz theorem), see \cite[p. 531]{Tadmor} and \cite{DMComptes}. Note that Tadmor actually used the notation  $V^{p, q} (\log V)^\alpha$ rather than $T^{p, q} \log^\alpha$, however we prefer to use the latter to emphasize the connection between Morrey and Tadmor spaces.} scale of function spaces $T^{\log, \alpha}$, which is obtained when the $\ell^\infty$-norm on the index set $I$ in \eqref{MorreyEquiv} is replaced by the bigger $\ell^2$-norm, specifically, 
$$
	\|\omega\|_{T^{\log, \alpha}} = \sup_{\Pi} \left\{\sum_{i \in I} \bigg[ (1 - (\log |Q_i|)_-)^{\alpha} |\omega|(Q_i) \bigg]^2  \right\}^{\frac{1}{2}}. 
$$
Equipped now with Tadmor spaces, the log-regularity $1/2$ can be achieved  as a threshold for the lack of concentrations in Open Problem \ref{OPDM}. To be more precise, the following result may be found in \cite[Corollary 4.1]{Tadmor}.

\begin{thm}\label{THMT} 
Assume that $\{u^\varepsilon\}$ is an approximation family of the $2$D Euler equations \eqref{Euler} such that the  corresponding family of vorticities
\begin{equation*}
 \{\omega^\varepsilon\} \quad  \text{is uniformly bounded in} \quad L^\infty ([0, T]; T^{\log, \alpha}) \quad \text{for some} \quad \alpha > \frac{1}{2}.
\end{equation*}
Then $\{\omega^\varepsilon\}$ defines a bounded family of vortex sheets  \eqref{Vortex} and
\begin{equation*}
\lim_{\varepsilon \to 0} u^\varepsilon = u \quad \text{in} \quad L^\infty ([0, T]; L^2_{loc}),
\end{equation*}
where $u$ is a weak solution to \eqref{Euler}. 
\end{thm}

The previous result consists of a remarkable extension of Theorem \ref{ThmDMConser} to the expected full range of regularity for strong convergence $\alpha > 1/2$. On the other hand, since the following trivial embedding holds 
$$
	T^{ \log, \alpha} \hookrightarrow M^{\log, \alpha},
$$
 regularity assumption from Theorem \ref{THMT}  seems to be slightly stronger than classical Morrey regularity as considered in Theorem \ref{ThmDMConser}.  Despite the fact that Tadmor's result does not apparently resolve Open Problem \ref{OPDM}, which was already pointed out by the author in \cite[p. 519 and the discussion after eq. (3.5)]{Tadmor} and again in \cite{TadmorLecture}, it gives another convincing reason (in addition to  the above mentioned result of Majda for $\mathcal{M}_+ \cap M^{\log, 1/2}$) to support the belief that $M^{\log, 1/2}$ is the critical regularity space for concentrations to occur.

 An interesting feature somehow implicit in Tadmor's method  is that special sets  (packings) of the underlying domain have something to say about   lack of concentrations. This observation served as an impetus to the first author and Milman \cite{DM} to develop a sparse approach to $H^{-1}$-stability (see Section \ref{Sstab}.) 
 
 \subsubsection{\emph{\textbf{The sparse method of Dom\'inguez and Milman: Energy conservation is encoded in sparse cubes.}}}\label{SecSpar}
It is well known from the pioneering work of Calder\'on and Zygmund in the 50's that packings play a central role in a variety of questions arising in classical harmonic analysis. On the other hand, modern developments within the last ten years in harmonic analysis abandon packings and rely instead on sparse families of cubes. Loosely speaking, sparse cubes are not necessarily packings, but possible overlappings may be controlled in a sharp fashion.  This principle (i.e., going from packings to sparse) has shown to have powerful consequences beyond the classical Calder\'on--Zygmund theory, but primarily focused on sharp estimates for operators via the so-called sparse domination, see \cite{LN}. Very recently, the first author and Milman \cite{DM} proposed to incorporate sparseness as a tool to characterize energy conservation  in Euler equations. Since sparseness will become a key ingredient in our later arguments, it would be convenient to make a quick review on the main achievements from  \cite{DM}. Next we recall the definition of sparse cubes. 

\begin{defn}[Sparse cubes]\label{DefSp}
	\normalfont A dyadic family of cubes $(Q_i)_{i \in I}$ is said to be \emph{sparse} if for every $Q_i$ there exists a measurable set $E_{Q_i} \subset Q_i$ such that $(E_{Q_i})_{i \in I}$ are pairwise disjoint and\footnote{For the purposes of this paper, the constant $1/2$ is unessential and may be replaced by another fixed $c \in (0, 1)$, say $c|Q_i| \leq |E_{Q_i}|$.} $|Q_i|/2 \leq |E_{Q_i}|$. We denote by $S$ the set of all sparse families of cubes.
\end{defn}

 Clearly $\Pi \subset S$, but $S$ contains families of overlapping cubes. To fix ideas,  a model example in $1$D of sparse intervals that are not a packing is given by dyadic intervals starting at the origin $\{[0, 2^{-k}]: k \geq 0\} \in S \backslash \Pi$.   
 
 The main goal of \cite{DM} is to show that the $H^{-1}$-stability method of Lopes Filho, Nussenzveig Lopes and Tadmor can be sharpened in terms of the so-called \emph{sparse indices} defined by 
 \begin{equation}\label{SI1}
 	s_n(\omega) =  \sup_{(Q_i)_{i \in I} \in S} \bigg(\sum_{i \in I: \ell(Q_i) < 2^{-n}} |\omega|(Q_i)^2  \bigg)^{\frac{1}{2}}, \qquad n \in \N,
 \end{equation}
 if $\omega \in \mathcal{M}$. Then the following numbers are introduced in a natural way
 \begin{equation}\label{SIn4}
 	s_n(\{\omega^\varepsilon\}) = \sup_{\varepsilon > 0}  s_n(\omega^\varepsilon)
 \end{equation}
 and
 \begin{equation}\label{SI2}
 	     s_n(X) = \sup_{\|\omega\|_X < 1} s_n(\omega)
 \end{equation}
 if $X \subset \mathcal{M}$ is a normed space. 
 In short, the indices $s_n(\{\omega^\varepsilon\})$ and $s_n(X)$  \emph{measure the concentration on sparse cubes} of $\{\omega^\varepsilon\}$ and the elements in the unit ball of $X$, respectively. Armed now with these  indices, Theorem \ref{THMLNT} admits the following improvement via sparseness  (cf. \cite[Corollary 1]{DM}). 
 
\begin{thm}\label{ThmDM2}
Assume that $\{u^\varepsilon\}$ is a family of approximate solutions of the $2$D Euler equations \eqref{Euler} such that 
\begin{equation}\label{CSp}
	\lim_{n \to \infty} s_n(\{\omega^\varepsilon\}) =0. 
\end{equation}
Then 
\begin{equation*}
\lim_{\varepsilon \to 0} u^\varepsilon = u \qquad \text{in} \qquad L^\infty ([0, T]; L^2_{loc}),
\end{equation*}
 where $u$ is a weak solution to \eqref{Euler}.	In addition, if $u$ is a physically realizable weak solution  then $u$ conserves the energy \eqref{Conservative}. 
\end{thm} 

This result contains precious information on the lack of concentrations of approximate solutions of Euler equations since \eqref{CSp} reveals us that energy preservation of Euler solutions is specially relevant in sparse families of cubes, which enjoy a rather simple geometric structure. In an informal manner, when compared with Tadmor's approach (see Section \ref{SectionT}), Theorem \ref{ThmDM2} claims that sparse cubes make a better job than packings! In practice, computability of $s_n(X)$ (and $s_n(\omega)$) is easy to implement for all relevant choices of $X$ (cf. \cite[Section 5]{DM}; see also \eqref{EstimM} and \eqref{N2662} below). As a byproduct,  Theorem \ref{ThmDM2} does not only recover in a unifying fashion all previously known existence results without concentration stated in \cite{DiPernaMajda, LNT, Tadmor}, but it also provides with a new constructive approach to regularity classes when combined with the  extrapolation theory \cite{JM} of Jawerth and Milman  (cf. \cite{DM} for further details.) 

Despite the fact that \cite{DM} does not give a definite answer to Open Problem \ref{OPDM}, the following fact underlying its main results (in particular, Theorem \ref{ThmDM2}) will play a fundamental role in the arguments of this paper: the asymptotic decay of $s_n(X)$ gives quantitative information on the degree of $H^{-1}$-compactness.  In particular, sparse indices of Morrey and Tadmor spaces (cf. \eqref{SI2}) can be estimated by 
 \begin{equation}\label{EstimM}
 	s_n(M^{\log, \alpha}) \leq c n^{\frac{1-\alpha}{2}} \qquad \text{if} \qquad \alpha > 1
 \end{equation}
 and
 \begin{equation}\label{N2662}
  s_n(T^{ \log, \alpha}) \leq c n^{\frac{1}{2} - \alpha} \qquad \text{if} \qquad \alpha > \frac{1}{2}. 
  \end{equation}
  These estimates together with the decay to zero imposed in \eqref{CSp} give  convincing explanations behind  regularity assumptions required in Theorems \ref{ThmDMConser} and \ref{THMT}  (i.e.,  $\alpha > 1$ and $\alpha > 1/2$ respectively).  More importantly, in striking contrast with classical Majda's result for $\mathcal{M}_+ \cap M^{\log, 1/2}$ and Theorem \ref{THMT}, \eqref{EstimM} is the first (but  still inconclusive) evidence that $\alpha > 1$ may be the optimal regularity for concentrations in Open Problem \ref{OPDM}. This conjecture is indeed confirmed by Theorems \ref{ThmMain} and \ref{ThmMain2}, see Section \ref{sparse_index} below.

\subsubsection{\emph{\textbf{Further related results.}}} 
We briefly mention that another interesting approach to $H^{-1}$-stability has been recently proposed by Lanthaler, Mishra, and Par\'es-Pulido \cite{Lan} in terms of the so-called structure functions\footnote{We observe that structure functions from \cite{Lan} are equivalent to classical moduli of smoothness in $L^2$.} of $\{u^\varepsilon\}$. In particular, this tool is employed by Lanthaler in \cite[Theorem 2.13]{Lan2} to derive another proof of Theorem \ref{ThmDMConser} via sharp estimates between structure functions and the maximal vorticity function. However, the method breaks down if $\alpha \leq 1$ since the estimate stated in \cite[Corollary 2.14]{Lan2} essentially requires $\alpha > 1$. Furthermore, several numerical simulations for lack of concentrations in  dynamics of unsigned vortex sheets are performed in \cite{Lan2}.  

On the other hand, Jiu and Xin \cite{JX, JX2} thoroughly investigated strong convergence of approximate solutions, generated by both smoothing the initial data and viscous approximations, to $3$D axisymmetric Euler equations with vortex sheets initial data. Loosely speaking, the authors show that if concentrations arise in the limit process then these must happen in a region outside the symmetry axis.  Among other results, the counterparts of Theorem \ref{ThmDMConser} involving vortex sheets under $\log$-Morrey regularity\footnote{Recall that  critical regularity for vortex sheets in $3$D is measured by the following variant of \eqref{Morrey}:
$
	\sup_Q \frac{(1-(\log |Q|)_-)^\alpha}{|Q|} \, |\omega|(Q).
$
See \cite{LNT} and \cite{DM}. 
} of order  $\alpha > 1$ and its optimality assertion for signed vortex sheets with  $\log$-Morrey regularity $\alpha=1/2$ are established in \cite[Section 3 and 4]{JX}.  As mentioned in \cite[p. 388]{JX} and \cite[p. 49]{JX2}, the gap between these two results,  i.e.,  the analog of Open Problem \ref{OPDM} for $3$D axisymmetric vortex sheets, remains as an open question. We believe that it is worth investigating whether methodology of Theorems \ref{ThmMain} and \ref{ThmMain2} proposed in this paper can be adequately adapted to resolve also the $3$D gap problem for axisymmetric vortex sheets as raised in \cite{JX, JX2}.

\subsection{Details about the proofs of Theorems \ref{ThmMain} and \ref{ThmMain2}} \hfill\\

\subsubsection{\textbf{\emph{Some challenges in the proof}}}
The original method of Theorem \ref{ThmDMConser} given in \cite{DiPernaMajda} to show lack of concentrations in $M^{\log, \alpha}, \, \alpha > 1$,  relies on regularity theory for elliptic PDE, namely,  H\"older-type estimates for the  streamfunction $\psi$ satisfying $\Delta \psi = \omega$ and $u= \nabla^\perp \psi$. On the other hand, the counterexample from \cite[Proposition 3.1]{DiPernaMajda1*} (see also \cite[pages 310 and 311]{DiPernaMajda} for a summary) exhibiting concentrations in $M^{\log, 1/2}$ is defined as the scalings:
\begin{equation}\label{Scaling}
	u^\varepsilon (x) = \bigg(\log \bigg(\frac{1}{\varepsilon} \bigg) \bigg)^{-1/2} \varepsilon^{-1} u_0 \bigg(\frac{x}{\varepsilon} \bigg),
\end{equation}
where $u_0$ is the steady rotating eddy induced by a smooth radial vorticity $\omega_0$ with compact support  and non-zero net circulation, i.e., 
\begin{equation}\label{Net}
	\int_0^\infty s \omega_0(s) \, ds \neq 0. 
\end{equation}
It is not hard to see that the scaling prefactor $(\log (1/\varepsilon))^{-\alpha}$ with $\alpha = 1/2$ in \eqref{Scaling} is the only possible choice of $\alpha$ that guarantees simultaneously, for corresponding vorticities $\omega^\varepsilon = \text{curl } u^\varepsilon$, 
$$
	\omega^\varepsilon \in M^{\log, \alpha} \qquad  \text{and} \qquad  \omega^\varepsilon \in  H^{-1}_{loc}
$$ 
uniformly in $\varepsilon > 0$. Accordingly,  velocity fields  \eqref{Scaling} cannot work beyond the setting of $M^{\log, 1/2}$ proposed in \cite{DiPernaMajda1*}. 

Another primary obstacle to be overcome comes from the assumption \eqref{Net}. That is, in light of the results from \cite{M93},  trying to show concentrations beyond $M^{\log, 1/2}$ requires to work with mixed-sign vorticities and then \eqref{Net} seems to be incoherent since zero net circulations should be expected in this setting.

\subsubsection{\textbf{\emph{Strategy of the proof}}} 
The main idea underlying the proofs of both Theorems \ref{ThmMain} and \ref{ThmMain2} is a novel approach that enables  to generate wild behaviour of solutions of $2$D Euler equations on adequate geometrical decompositions of their corresponding domains in terms of sparse cubes (see Definition \ref{DefSp}). To make this assertion more precise, we recall that on account of the recent method introduced by the first author and Milman \cite{DM} (cf. Section \ref{SecSpar}), in order to  characterize  energy conservation/anomalous dissipation for solutions to Euler equations is enough to control their behaviour within the distinguished geometry given by sparse cubes.  Motivated by this sparse principle, we are going to construct a  Cantor set 
\begin{equation}\label{4n}
E = \bigcap_{k=0}^\infty E_k
\end{equation}
within the unit square, whose $k$-th generation $E_k$ is formed by $2^{2 k}$ dyadic cubes with side length $2^{-2^{2 k}}$. In particular
\begin{equation}\label{IntroMor3}
E_{k+1} \subset E_k.
\end{equation}
 Consider now the family consisting of all these cubes, namely, 
\begin{equation}\label{PSp1}
	\mathcal{Q} = \{Q : Q \in E_k, \,   k \geq 0\}. 
\end{equation}
It is obvious that $\mathcal{Q}$ does not define a packing $\mathcal{Q} \not \in \Pi$ (cf. \eqref{IntroMor3}). However, we will show  that $\mathcal{Q}$ is a sparse family of cubes, that is,  possible overlappings of their elements $Q \in \mathcal{Q}$ can be controlled in a sharp way. 

Once the desired sparse family $\mathcal{Q}$ has been settled (cf. \eqref{PSp1}), the proof of Theorem \ref{ThmMain} seeks for a measure $\omega$ having adequate decay rates for circulation on each $Q \in \mathcal{Q}$ leading to $\omega \not \in H^{-1}_{loc}$ (i.e., $\omega$ is not a vortex sheet).  To this end, we define the sequence $\{\omega^k\}$ as an equidistributed number of Dirac masses $\delta_{x_{k, m}}$ related to the centers $x_{k, m}$ of each $Q \in \mathcal{Q} \cap E_k$, more precisely, 
\begin{equation}\label{IntroMor2}
	\omega_k = \sum_{m=1}^{2^{2 k}} 2^{-2 k} \delta_{x_{k, m}}. 
\end{equation}
Observe that $\omega_k$ is a probability measure on the unit square and, in addition, 
 elementary computations show that
\begin{equation}\label{34n}
	\omega_k(Q) = c (-\log |Q|)^{-1} \qquad \text{for every} \qquad  Q \in E_k.
\end{equation}
 Let us then define   $\omega$ to be the weak-limit of the sequence $\{\omega_k\}$. As a consequence of \eqref{34n}, we will achieve the desired regularity $\omega \in M^{\log, 1}$, despite the fact that\footnote{Observe that, by \eqref{IntroMor2},  $\omega_k (Q) = 2^{-2 k}$  for every  cube $Q$ such that $x_{k, m} \in Q$ with side length sufficiently small. Then $\omega_k \not \in M^{\log, 1}$ since $\lim_{|Q| \to 0} \omega_k(Q) = 2^{-2 k} \neq 0$.  If one is so inclined, replacing Dirac masses with suitable scalings of the Lebesgue measure gives a sequence $\{\tilde{\omega}_k\}_{k \in \mathbb{N}}$ with the same weak-star limit $\omega$ and such that $\tilde{\omega}_k \in M^{\log, 1}$.} $\omega_k \not \in M^{\log, 1}$. 
 
 On the other hand, the sparseness of the underlying family of cubes $\mathcal{Q}$ will play a key role in the proof of  $\omega \not \in H^{-1}_{loc}$. Indeed, given $Q \in E_k$, let us denote by  $\widehat{Q}$ the unique ancestor of $Q$ in the previous Cantor generation  $E_{k-1}$. Then we are going to construct an increasing tower of sparse cubes going from $Q$ to $\widehat{Q}$ that enables us to compute sparse indices of $\omega$ (cf. \eqref{SI1}) in a relatively simple way. In fact, this sparse construction turns out to be optimal, in the sense that it yields  worst possible behaviour of sparse indices  $s_n(\omega) = \infty$. In particular, this argument concludes  $\omega \not \in H^{-1}_{loc}$.

Concerning the proof of Theorem \ref{ThmMain2}, similar ideas as above can still be implemented under some adequate modifications. To be more precise, the main difference with respect to Theorem \ref{ThmMain} comes from the fact that the sequence $\{\omega_k\}$ is now required to be formed by vortex sheets. This makes computations (but not arguments themselves) involved in the proof of Theorem \ref{ThmMain2} to be slightly more technical than those corresponding to Theorem \ref{ThmMain}. At this stage, we are inspired by the nice construction of vortex sheets due to Greengard and Thomann \cite{GT} in connection with convergence issues of the weak$^*$ defect measure $u^\varepsilon \otimes u^\varepsilon$ as $\varepsilon \to 0$. Then we propose to replace the role of the Dirac masses $2^{-2 k} \delta_{x_{k, m}}$ in the definition of $\omega_k$ given in  \eqref{IntroMor2} by scalings adapted to each cube $Q \in E_k$ of a fixed radial vortex patch $\omega_0$  with zero net circulation
\begin{equation}\label{IntroMor4}
	\int_0^\infty s \omega_0(s) \, ds=0. 
\end{equation}
This is in sharp contrast to the counterexample \eqref{Scaling} of DiPerna and Majda working with $M^{\log, 1/2}$ for which   \eqref{Net} is satisfied. 
From the classical Biot-Savart law (i.e., $u = k \ast \omega$ where $k$ is the Biot-Savart kernel), one can compute explicitly the exact steady solution $u_0$ to \eqref{Euler}  with  $\text{curl } u_0 = \omega_0$. In particular, \eqref{IntroMor4} guarantees that the support of $u_0$ is contained in the support of the given vortex patch $\omega_0$.  Having arrived at this point,  we can construct the sequence $\{u_k\}$ with $\text{curl } u_k = \omega_k$ such that every $u_k$ is formed by $2^{2 k}$ non intersecting  copies of the steady exact solution $u_0$ adequately scaled to each cube $Q$ in the $k$-th Cantor generation $E_k$ of $E$, see \eqref{4n}. To prove that $\{u_k\}$ and $\{\omega_k\}$ fulfil the list of properties stated in Theorem \ref{ThmMain2}, we will again make an essential use of the sparse  geometry of the cubes $Q \in E_k$, the supports of vortex patches $\omega_k$.

As a remark, we mention that the process of gluing vortex patches proposed in the above Cantor-type construction is essential in order to achieve  claims given in Theorem \ref{ThmMain2}. In particular, this construction produces a completely different outcome when compared with the one of Greengard and Thomann \cite{GT}: Scalings from \cite{GT} are adapted to the family of \emph{all} dyadic cubes contained in the unit square and, as a consequence, $\theta$ coincides with the Lebesgue measure on unit square, showing the worst possible behaviour of $\theta$ in the sense that all measurable sets with positive measure are concentration sets. This is not the case in our sparse construction, which only induces wild behaviours of Euler solutions on properly selected sparse families of cubes and, in particular, $\theta$ coincides with the  measure $\omega$ of logarithmic type obtained in Theorem \ref{ThmMain}.  As already stressed above, the underlying sparse structure  plays a fundamental role in the method of proof of Theorem  \ref{ThmMain2}.   

\subsection{Interpretation of Theorems \ref{ThmMain} and \ref{ThmMain2} in terms of $H^{-1}$-stability} \hfill\\

As recalled in Subsection \ref{Sstab}, to solve Open Problem \ref{OPDM} can be reduced to study  optimal forms of the following compactness assertion
\begin{equation}\label{CMorreyNew}
	M^{\log, \alpha}  \overset{compactly}{\hookrightarrow} H^{-1}_{loc}, 
\end{equation}
as well  as its variants for special $M^{\log, \alpha}$-sequences of approximate solutions. 
Up to now, the best available results were obtained in \cite{LNT} and \cite{M93} showing  the validity of \eqref{CMorreyNew} under $\alpha > 1$ and its failure with $\alpha = 1/2$, respectively.  To close the gap between these two results has remained as an open question, which has been explicitly raised in the works of Lopes Filho, Nussenzveig Lopes, and Tadmor \cite{LNT} and Tadmor \cite{Tadmor, TadmorLecture}. 

Results from this paper gives a final answer to \eqref{CMorreyNew}: 
$$
	M^{\log, \alpha}  \overset{compactly}{\hookrightarrow} H^{-1}_{loc}  \iff \alpha > 1. 
$$
Indeed, Theorem \ref{ThmMain} provides with the  much stronger statement:
\begin{equation}\label{Dfv}
	M^{\log, 1}  \not \hookrightarrow H^{-1}_{loc}. 
\end{equation}
On the other hand, Theorem \ref{ThmMain2} asserts that \eqref{Dfv} does not admit an improvement even when attention is restricted to the special  $M^{\log, 1}$-set formed by classical solutions of $2$D Euler equations. In particular, 
\begin{equation}\label{Dfv2}
	M^{\log, 1} \cap H^{-1}_{loc}  \overset{compactly}{\not\hookrightarrow} H^{-1}_{loc}. 
\end{equation}

\subsection{Organization of the paper} \hfill\\

The proofs of Theorems \ref{ThmMain} and \ref{ThmMain2} are given in Sections \ref{Section3} and \ref{Section4}, respectively.

%%%%%%%%%%%%%
\section{Proof of Theorem \ref{ThmMain}}\label{Section3}

\subsection{Construction of the Cantor set}\label{Step1} Consider the sequence $\{l_k\}_{k \geq 0}$ defined by
\begin{equation}\label{Sequence}
	l_k = 2^{-2^{2 k}}, \quad k \in \N, \qquad l_0 = 1.
\end{equation}
Note that $2 l_{k+1} < l_k$ and, in particular, $\{l_k\}_{k \geq 0}$ is decreasing. 
Let $I_0 = [0, 1]$ and define $I_1$ as the subset of $I_0$ which is obtained by removing an open interval of length $l_0 - 2 l_1$ in the middle of $I_0$, more precisely, $$I_1 = [0, l_1] \cup [1-l_1, 1].$$ Then $I_1$ is formed by $2$ dyadic intervals of length $l_1$. We apply the same procedure to each of these intervals to define $I_2$, that is, $$I_2 = [0, l_2] \cup [l_1-l_2, l_1] \cup [1-l_1, 1-l_1 + l_2] \cup [1-l_2, 1].$$
Then $I_2$ is formed by $2^2$ dyadic intervals of length $l_2$. Iterating the construction, one can obtain a sequence $\{I_k\}_{k \geq 0}$ of subsets of $I_0$ such that $I_{k+1} \subset I_k$, where each $I_k$ is formed by $2^k$ dyadic intervals of  length $l_k$. 

The previous $1$D construction can be extended to $2$D in a natural way: Let $E_k$ be the subset of $Q_0 = [0, 1]^2 = I_0 \times I_0$ which is formed by the Cartesian products of the $I_k$'s. Then
\begin{equation}\label{nested}
E_{k+1} \subset E_k
\end{equation}
 and each $E_k$ is the union of $2^{2 k}$ cubes $Q_{2^{2 k}, m}$ within the dyadic generation $2^{2 k}$ of the original cube $Q_0$, namely, 
\begin{equation}\label{Es}
	E_k = \bigcup_{m=1}^{2^{2 k}} Q_{2^{2 k}, m}
\end{equation}
with $\ell(Q_{2^{2 k}, m}) = l_k =2^{-2^{2 k}}$. 

The  Cantor set $E$  relative to $\{E_k\}_{k \geq 0}$ is defined by
\begin{equation}\label{Cantor}
	E = \bigcap_{k =0}^\infty E_k. 
\end{equation}

\begin{figure}[H]
\begin{center}
\begin{tikzpicture}

\draw[very thick] (0, 0) rectangle (10,10);

\draw[very thick] (0, 0) rectangle (3,3);
\draw[very thick] (7, 0) rectangle (10, 3);
\draw[very thick] (0, 7) rectangle (3, 10);
\draw[very thick] (7, 7) rectangle (10, 10);

\draw[very thick] (0, 0) rectangle (4/5, 4/5);
\draw[very thick] (3-4/5, 0) rectangle (3, 4/5);
\draw[very thick] (0, 3-4/5) rectangle (4/5, 3);
\draw[very thick] (3-4/5, 3-4/5) rectangle (3, 3);

\draw[very thick] (7, 0) rectangle (7 + 4/5, 4/5);
\draw[very thick] (10-4/5, 0) rectangle (10, 4/5);
\draw[very thick] (7, 3-4/5) rectangle (7+4/5, 3);
\draw[very thick] (10-4/5, 3-4/5) rectangle (10, 3);

\draw[very thick] (0, 10-4/5) rectangle (4/5, 10);
\draw[very thick] (0, 7) rectangle (4/5, 7 + 4/5);
\draw[very thick] (3-4/5, 7) rectangle (3, 7 + 4/5);
\draw[very thick] (3-4/5, 10-4/5) rectangle (3, 10);

\draw[very thick] (7, 7) rectangle (7 + 4/5, 7 + 4/5);
\draw[very thick] (10-4/5, 7) rectangle (10, 7 + 4/5);
\draw[very thick] (7, 10-4/5) rectangle (7 + 4/5, 10);
\draw[very thick] (10-4/5, 10-4/5) rectangle (10, 10);
\end{tikzpicture}
\end{center}
\caption{Construction of the Cantor set $E$.}
\end{figure}
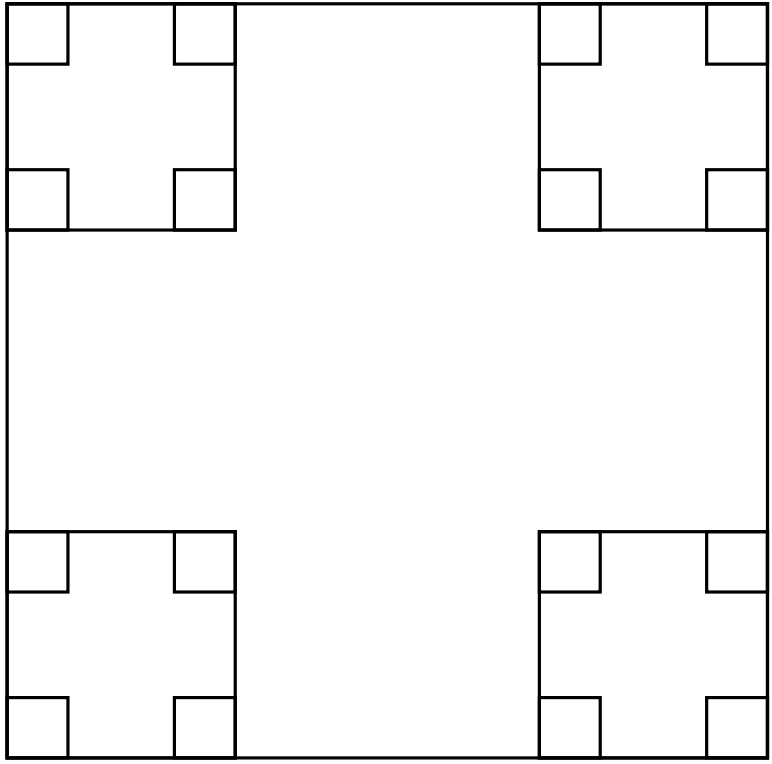

\subsection{Construction of the measure associated with the Cantor set}\label{Step2} We wish to construct $\omega \in \mathcal{M}^+$ relative to the Cantor set $E$ given by \eqref{Cantor} with $\omega(E)=1$. This can be done via standard limiting arguments. Indeed, for every $k$ we define the measure $\omega_k$ as
\begin{equation}\label{Measure}
	\omega_k = \sum_{m=1}^{2^{2 k}} 2^{-2 k} \delta_{x_{k, m}}, 
\end{equation}
where $\delta_{x_{k, m}}$ denotes the Dirac mass related to $x_{k, m}$, the center of the cube $Q_{2^{2 k}, m}$ given in  \eqref{Es}. In particular, $\omega_k$ is a probability measure on $E_k$ since
$$
	\omega_k (E_k) = \sum_{m=1}^{2^{2 k}} \omega_k (Q_{2^{2 k}, m}) = 2^{-2 k} 2^{2 k} = 1. 
$$
As a consequence of the Riesz representation theorem for the dual of space of continuous functions, $\{\omega_k\}_{k \geq 0}$ converges  weakly (possibly passing to a subsequence) to a measure $\omega \in \mathcal{M}^+$.  

Let $k \geq 0$ and $m \in \{1, \ldots, 2^{2 k}\}$. By construction, it is not hard to check that  
$$
	\omega_j(Q_{2^{2 k}, m}) = \omega_k (Q_{2^{2 k}, m})  \qquad \text{for} \qquad k \leq j
$$
and, in particular, 
\begin{equation}\label{P1}
	\omega(Q_{2^{ 2k}, m}) = \lim_{j \to \infty} \omega_j(Q_{2^{ 2k}, m}) =  \omega_k (Q_{2^{2 k}, m}) = 2^{-2 k},
\end{equation}
where we have used  \eqref{Measure} in the last step. 
Furthermore,  $\omega$ is a probability measure on the Cantor set $E$ (cf. \eqref{Cantor}). Indeed, this follows from \eqref{P1} since
$$
	\omega(E) = \lim_{k \to \infty} \omega (E_k) = \lim_{k \to \infty} \sum_{m=1}^{2^{2 k}} \omega(Q_{2^{ 2k}, m}) = \lim_{k \to \infty} 2^{2 k} 2^{-2 k} = 1. 
$$
This gives the desired measure $\omega$ on the Cantor set $E$. 

\subsection{$\omega$ is a Morrey measure in $M^{\log, 1}$} Let $Q$ be a (not necessarily dyadic) cube in $Q_0$ with side-length $\ell(Q)$. Then there exists a unique $j = j_{Q}$ such that 
\begin{equation}\label{P2}
l_{j+1} \leq \ell(Q) < l_j.
\end{equation}
 This is possible since $\{l_k\}_{k \geq 0}$ is a decreasing sequence with $\lim_{k \to \infty} l_k = 0$, cf. \eqref{Sequence}.  A simple application of the pigeon-hole principle yields that there exists an absolute constant $C \geq 1$ such that $Q$ intersects at most $C$ cubes $Q_{2^{j}, m}$ from the collection $E_j$. Since $\omega$ is uniform on all cubes contained in $E_j$ (cf. \eqref{P1}), we arrive at\footnote{Given two non-negative quantities $A$ and $B$, the standard notation $A \lesssim B$ means that there exists a constant $c$, independent of all essential parameters, such that $A \leq c B$. We write $A \approx B$ if $A \lesssim B \lesssim A$.} 
\begin{equation}\label{P3}
	\omega(Q) \leq C \omega (Q_{2^{j}, m})  = C 2^{-2 j} \lesssim  (1-\log |Q|)^{-1},
\end{equation}
where the last estimate follows from the fact that $1-\log |Q| \approx 2^{2 j}$ (recall \eqref{P2}). 
The fact that the bound \eqref{P3} is uniform with respect to all cubes $Q \subset Q_0$ enables to conclude that  $\omega \in M^{\log, 1}$. 

\subsection{Geometrical construction of sparse sets}\label{Section4n} Our next goal is to construct an adequate family of sparse cubes that is able to capture the energy concentration of the measure $\omega$ introduced in Subsection \ref{Step2}. To proceed with, we are going  to augment in a clever way   the family of the dyadic cubes involved in the construction of the Cantor set $E$ given in Subsection \ref{Step1}. Recall that, for every $k \in \N$, the set $E_k$ is formed by $2^{2 k}$ dyadic cubes $Q_{2^{2 k}, m}$ with $\ell(Q_{2^{2 k}, m}) = 2^{-2^{2 k}}$, cf. \eqref{Es}. In particular, these cubes are pairwise disjoint. Without loss of generality, we shall only argue for one of these cubes, $Q_{2^{2 k}, 1}$,  but the arguments can be carried out to all cubes $Q_{2^{2 k}, m}$ for $m= 1, \ldots, 2^{2 k}$. By construction, $Q_{2^{2 k}, 1}$ contains exactly 4 cubes of the next generation $E_{k+1}$ involved in the definition of the Cantor set $E$. Choose only one of these cubes, which is denoted by $Q_{2^{2 (k+1)}, 1}$. To simplify the exposition, we may think that $Q_{2^{2(k+1)}, 1}$ and its ancestor $Q_{2^{2 k}, 1}$ have common lower left vertex.  In particular,   we have $$Q_{2^{2(k+1)}, 1} \subset Q_{2^{2 k}, 1} \cup E_{k+1}.$$
Let us denote by
$$
	\{Q_{j, 1} : j = 2^{2 k}, \ldots, 2^{2(k+1)}\}
$$
the full tree of dyadic cubes with common lower left vertex that expands from $Q_{2^{2(k+1)}, 1}$ to $Q_{2^{2 k}, 1}$ and we collect all these cubes (except the original ancestor $Q_{2^{2 k}, 1}$) in the family $\mathscr{S}_{k, 1}$, namely,  
\begin{equation*}
\mathscr{S}_{k, 1} = \{Q_{j, 1} : j = 2^{2 k} + 1, \ldots, 2^{2(k+1)} \}. 
\end{equation*}
Note that
\begin{equation}\label{P5a}
	|\mathscr{S}_{k , 1}| = 2^{2 k}
\end{equation}
 and
\begin{equation}\label{P5}
	Q_{2^{2 k}, 1} \supset Q_{2^{2 k} +1, 1} \supset \cdots \supset Q_{2^{2(k+1)}, 1}. 
\end{equation}
In a similar fashion, the families of cubes $\mathscr{S}_{k, m}$ relative to $Q_{2^{2 k}, m}$ for any $m=1, \ldots, 2^{2 k}$  can be defined. Moreover, for every fixed $k$, the fact that the cubes $Q_{2^{2 k}, m}$ are pairwise disjoint yields that  the families $\mathscr{S}_{k, m}$ are also pairwise disjoint.  

Let 
$$
	\mathscr{S} = \bigcup_{k=0}^\infty \bigcup_{m=1}^{2^{2 k}}   \{Q : Q \in \mathscr{S}_{k,  m}\}.
$$
It is clear that $\mathscr{S} \not \in \Pi$ (in particular, the cubes in $\mathscr{S}_{k, 1}$ are not pairwise disjoint). However, it is easy to check that $\mathscr{S} \in S$. In fact, $\mathscr{S}$ may be viewed as the $2$D counterpart of the prototypical  sparse family $\{[0, 2^{-k}] : k \geq 0\}$ in the $1$D setting.

\begin{figure}[H]
\begin{center}
\begin{tikzpicture}
	
\draw[very thick] (0, 0) rectangle (10,10);	
	
	\draw[very thick] (0, 0) rectangle (2,2);
\draw[very thick] (8, 0) rectangle (10, 2);
\draw[very thick] (0, 8) rectangle (2, 10);
\draw[very thick] (8, 8) rectangle (10, 10);

\draw node at (5, 6) {$Q_{2^{2 k}, 1}$};

\draw[very thick, dashed] (0, 0) rectangle (5, 5);
%\draw[very thick, dashed] (5, 0) rectangle (10, 5);
%\draw[very thick, dashed] (0, 5) rectangle (5, 10);

\draw node at (1, 1) {$Q_{2^{2(k+1)}, 1}$};
%\draw node at (9, 1) {$Q_{2^{2(k+1)}, 2}$};
%\draw node at (1, 9) {$Q_{2^{2(k+1)}, 3}$};
%\draw node at (9, 9) {$Q_{2^{2(k+1)}, 4}$};

\draw node at (4,4.5) {$Q_{2^{2 k} + 1, 1}$};

\draw[very thick, dashed] (0, 0) rectangle (4, 4);
\draw node at (3.25,3.5) {$Q_{2^{2 k} + 2, 1}$};

\draw[very thick, dashed] (0, 0) rectangle (3, 3);
\draw node at (2.25,2.65) {$Q_{2^{2 k} + 3, 1}$};

\end{tikzpicture}
\end{center}
\caption{\label{Figure2}Construction of $\mathscr{S}$ in the $k$-th generation.}
\end{figure}
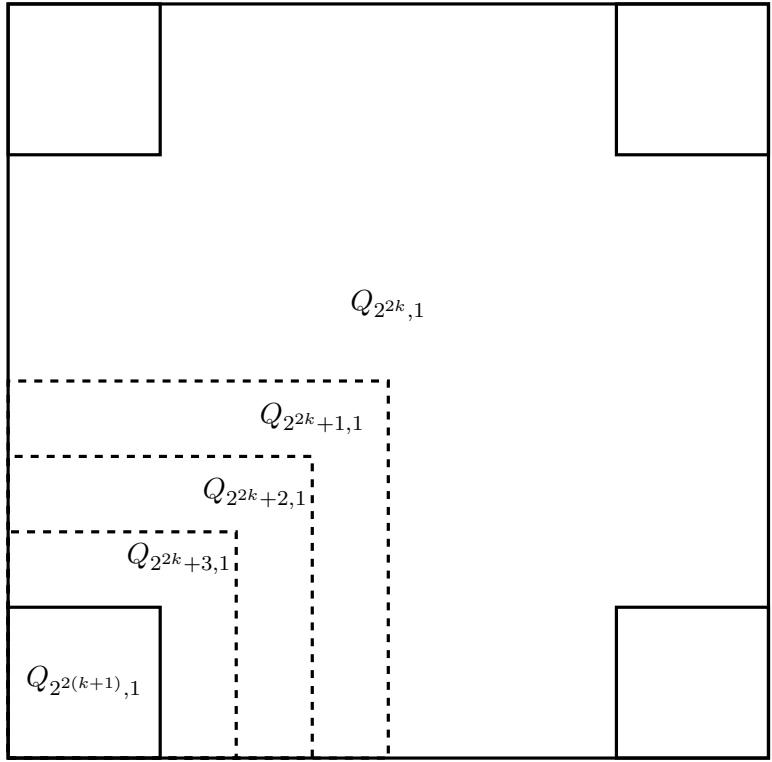

\subsection{Computability of sparse indices of $\omega$}\label{sparse_index}
In this section, we shall obtain sharp estimates for the sequence of sparse indices $s_n(\omega)$  relative to the measure $\omega$ introduced in Subsection \ref{Step2}. To do this, we will make use  of the sparse family  $\mathscr{S}$ given in Subsection \ref{Section4n}.

Recall that if $Q \in \mathscr{S}_{k, m} \subset \mathscr{S}, \, m \in \{1, \ldots, 2^{2 k}\},$ then $2^{-2^{2(k+1)}} \leq \ell(Q) < 2^{-2^{2 k}}$. For every $l \geq 0$, by  \eqref{SI1}, 
\begin{align}\label{P4}
	s_{2^{2 l}} (\omega)^2 &\geq \sum_{Q \in \mathscr{S}: \ell(Q) < 2^{-2^{2 l}}} \omega(Q)^2 \nonumber \\
	& = \sum_{k=l}^\infty  \sum_{Q \in \mathscr{S}: 2^{-2^{2(k+1)}} \leq \ell(Q) < 2^{-2^{2 k}}} \omega(Q)^2 \nonumber  \\
	&  =  \sum_{k=l}^\infty   \sum_{m=1}^{2^{2 k}}  \sum_{Q \in \mathscr{S}_{k, m}} \omega(Q)^2.
	\end{align}

Next we estimate $\omega(Q)$ for $Q \in \mathscr{S}_{k, m}$.  Again we make the non-restrictive assumption  $m=1$. It follows from \eqref{P5} that
$$
	Q_{2^{2(k+1)}, 1} \subset Q \subset Q_{2^{2 k}, 1},
$$
where $Q_{2^{2(k+1)}, 1} \in E_{k+1}$ and $Q_{2^{2 k}, 1} \in E_{k}$. In light of \eqref{P1}, 
$$
	\omega(Q_{2^{2 k}, 1} ) = 2^{-2 k}, \qquad  Q_{2^{2(k+1)}, 1} = 2^{-2(k+1)}.
$$
Then
\begin{equation}\label{P6}
	2^{-2 k} \approx \omega(Q_{2^{2(k+1)}, 1} ) \leq \omega(Q) \leq \omega(Q_{2^{2 k}, 1}) = 2^{-2 k}.
\end{equation}

From \eqref{P5a} and \eqref{P6}, we obtain 
$$
	\sum_{Q \in \mathscr{S}_{k, m}} \omega(Q)^2 \approx 2^{-4 k} 2^{2 k} = 2^{-2 k}
$$
for every $m \in \{1, \ldots, 2^{2 k}\}$. Inserting this estimate into \eqref{P4}, we derive
$$
	s_{2^{2 l}}(\omega)^2 \gtrsim \sum_{k=l}^\infty \sum_{m=1}^{2^{2 k}} 2^{-2 k} \approx \sum_{k=l}^\infty 2^{-2 k} 2^{2 k} = \sum_{k=l}^\infty 1 = \infty
$$
for all $l \geq 0$. 
As a consequence, noting that the sequence $s_n(\omega)$ is non-increasing, we conclude that 
$$
	s_n (\omega) = \infty,  \qquad \forall n \in \N. 
$$
\qed

\section{Proof of Theorem \ref{ThmMain2}}\label{Section4}

\subsection{Construction of vortex patches related to Cantor sets.} Recall the construction of the Cantor set $E$ given in Subsection \ref{Step1}. In particular, for every $k \in \N$, the set $E_k$ was defined in \eqref{Es} as the disjoint union of $2^{2 k}$ cubes $Q_{2^{2 k}, m}$ of side length $l_k = 2^{-2^{2 k}}$ (cf. \eqref{Sequence}). Denote by $x_{k, m}$ the center of $Q_{2^{2 k}, m}$.

For each $k \in \mathbb{N}$, we define  the following set of balls related to the Cantor sets $E_k$: 
$$
A_k = \{B(x_{k, m}, \delta_k) : m= 1, \ldots, 2^{2 k}\},
$$
 where the radii $\delta_k \in (0, 1)$ will be fixed later. We also introduce the set $B_k$ formed by all annuli  centered at $x_{k, m}$ with inner and outer radii $\sqrt{\delta_k}$ and $R_k$, respectively. To be more precise, 
 $$
 	B_k = \{\mathcal{C}_{k, m} : m = 1, \ldots, 2^{2 k}\},
 $$
 where
 $$
 	\mathcal{C}_{k, m} = B(x_{k, m}, R_k) \backslash B(x_{k, m}, \sqrt{\delta_k}). 
 $$
 In particular,  $R_k$ will be chosen small enough in order to guarantee that 
 \begin{equation}\label{Contain}
 	\mathcal{C}_{k, m} \subset B(x_{k, m}, R_k)  \subset Q_{2^{2 k}, m}. 
 \end{equation}

At the level $k$, we next construct the vorticity field $\omega_k$,  which is obtained when a suitable vortex patch   is copied  $2^{2 k}$ times according to the points $x_{k, m}$ generated in the construction of the Cantor set $E_k$. Specifically, we let
\begin{equation}\label{DefVort}
    \omega_k (x) = \Omega_k^+ \chi_{A_k}(x) + \Omega_k^- \chi_{B_k}(x),
\end{equation}
where the values $\Omega_k^+ > 0$ and $\Omega_k^- < 0$ will be specified later. The construction of $\omega_k$ is depicted in Figure \ref{Figure3} below. 

\vspace{1cm}

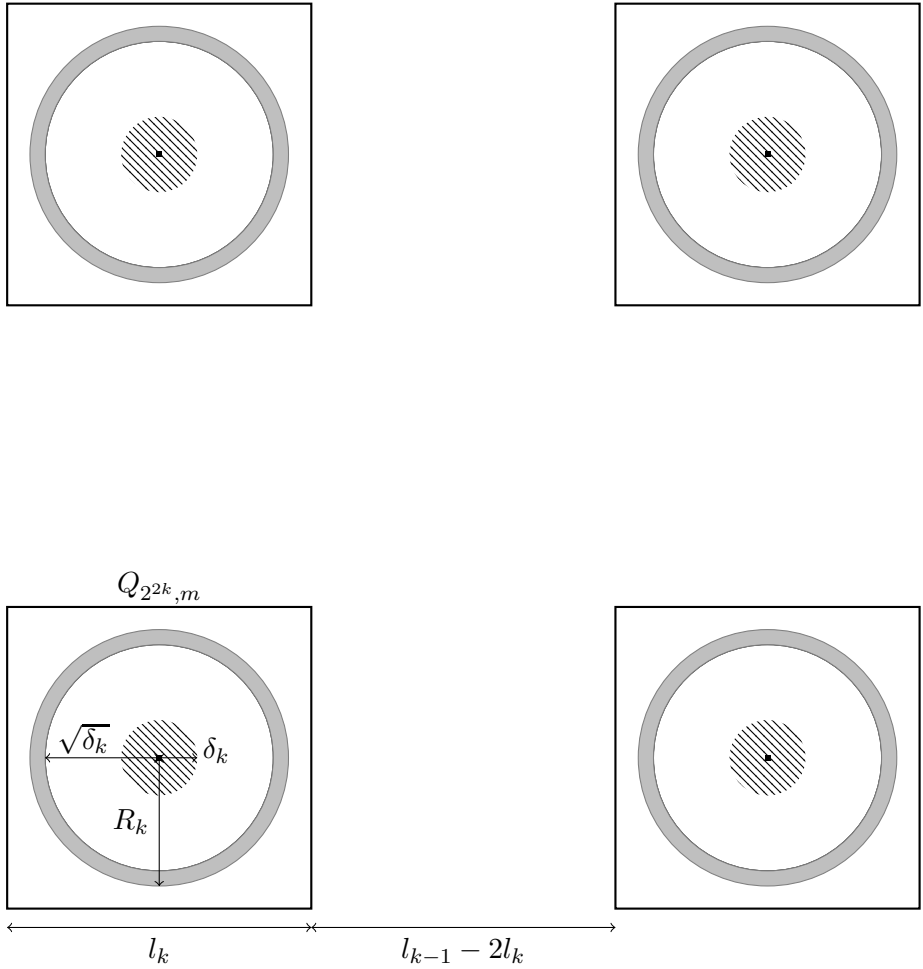
\begin{figure}[H]
\begin{center}
\begin{tikzpicture}
    % Define radius and center coordinates for the circles
    \def\Rk{1.5} % Outer radius
    \def\Rka{1.7}
    \def\rk{0.5} % Inner radius
    \def\pk{2}   % Distance from the corner to the center of the outer circle
    
    % Bottom-left circle
    \begin{scope}[shift={(0,0)}]
        \draw[thick] (0,0) rectangle (2*\pk, 2*\pk) node at (2*\pk/2, 2*\pk+0.25)  {$Q_{2^{2 k}, m}$}; % Square boundary
        \draw (1*\pk, 1*\pk) circle (\Rk); % Outer circle
    %    \draw[gray, fill = gray!50] (1*\pk, 1*\pk) circle (\rk); % Inner circle
        \path[pattern=north west lines,pattern color=black,even odd rule] (1*\pk, 1*\pk) circle (\rk);
        \draw[gray, fill=gray!50, even odd rule] (1*\pk, 1*\pk) circle (\Rk) circle (\Rka);
        \node[fill=black, inner sep=1pt] at (1*\pk, 1*\pk) {}; % Center dot
        % Radii and labels
        \draw[<->] (1*\pk, 1*\pk) -- (1*\pk+\rk, 1*\pk) node at (1*\pk + 0.75, 1*\pk + 0.1) {$\delta_k$};
        \draw[<->] (1*\pk, 1*\pk) -- (1*\pk-\Rk, 1*\pk) node at (1*\pk-1, 1*\pk+0.25)   {$\sqrt{\delta_k}$};
        \draw[<->] (0, -0.25) -- (2*\pk, -0.25) node[midway, below] {$l_k$};
        \draw[<->] (1*\pk, 1*\pk) -- (1*\pk, 1*\pk-\Rka) node[midway, left] {$R_k$};
        \draw[<->]  (2*\pk, -0.25) -- (8, -0.25) node[midway, below] {$l_{k-1}-2 l_k$};
    \end{scope}
    
    % Bottom-right circle
    \begin{scope}[shift={(8,0)}]
        \draw[thick] (0,0) rectangle (2*\pk, 2*\pk); % Square boundary
        \draw (1*\pk, 1*\pk) circle (\Rk); % Outer circle
                \path[pattern=north west lines,pattern color=black,even odd rule] (1*\pk, 1*\pk) circle (\rk);
      %  \draw[blue, fill = blue!50] (1*\pk, 1*\pk) circle (\rk); % Inner circle
                \draw[gray, fill=gray!50, even odd rule] (1*\pk, 1*\pk) circle (\Rk) circle (\Rka);
        \node[fill=black, inner sep=1pt] at (1*\pk, 1*\pk) {}; % Center dot
    \end{scope}
    
    % Top-left circle
    \begin{scope}[shift={(0,8)}]
        \draw[thick] (0,0) rectangle (2*\pk, 2*\pk); % Square boundary
        \draw (1*\pk, 1*\pk) circle (\Rk); % Outer circle
                \path[pattern=north west lines,pattern color=black,even odd rule] (1*\pk, 1*\pk) circle (\rk);
      %  \draw[blue, fill = blue!50] (1*\pk, 1*\pk) circle (\rk); % Inner circle
	 \draw[gray, fill=gray!50, even odd rule] (1*\pk, 1*\pk) circle (\Rk) circle (\Rka);
        \node[fill=black, inner sep=1pt] at (1*\pk, 1*\pk) {}; % Center dot
    \end{scope}
    
    % blueTop-right circle
    \begin{scope}[shift={(8,8)}]
        \draw[thick] (0,0) rectangle (2*\pk, 2*\pk); % Square boundary
        \draw (1*\pk, 1*\pk) circle (\Rk); % Outer circle
                \path[pattern=north west lines,pattern color=black,even odd rule] (1*\pk, 1*\pk) circle (\rk);
       % \draw[blue, fill = blue!50] (1*\pk, 1*\pk) circle (\rk); % Inner circle
	 \draw[gray, fill=gray!50, even odd rule] (1*\pk, 1*\pk) circle (\Rk) circle (\Rka);
        \node[fill=black, inner sep=1pt] at (1*\pk, 1*\pk) {}; % Center dot
    \end{scope}
\end{tikzpicture}

\end{center}
\caption{\label{Figure3}Construction of the vortex patches $\omega_k$ in the $k$-th generation of the Cantor set $E$.}
\end{figure}

\subsection{$\{\omega_k\}$ is uniformly bounded in $L^1$.}
For each $m= 1, \ldots, 2^{2 k}$, we have
\begin{align*}
	\int_{Q_{2^{2 k}, m}} \omega_k (x) \, dx &= \Omega_k^+ |B(x_{k, m}, \delta_k)| + \Omega_k^- |\mathcal{C}_{k, m}| \\
	& = \pi \Omega_k^+  \delta_k^2 + \pi \Omega_k^- (  R_k^2-\delta_k)
\end{align*}
Accordingly,  the choice of $\Omega_k^-$ given by
\begin{equation}\label{Omega_k-}
	\Omega_k^- = - \frac{\Omega_k^+ \delta_k^2}{R_k^2-\delta_k}
\end{equation}
leads that each patch of $\omega_k$ has zero total vorticity inside of its corresponding support, namely, 
\begin{equation}\label{ZeroVort}
	\int_{Q_{2^{2 k}, m}} \omega_k (x) \, dx  = 0. 
\end{equation}
In addition, as a consequence of \eqref{ZeroVort}, the total circulation of $\omega_k$ vanishes on $Q_0 = [0, 1]^2$:
$$
	\int_{Q_0} \omega_k (x) \, dx = \sum_{m=1}^{2^{2 k}} \int_{Q_{2^{2 k}, m}} \omega_k (x) \, dx =0. 
$$

Taking into account the value of $\Omega_k^-$ given in \eqref{Omega_k-},  we can compute the $L^1$ norm of $\omega_k$ as follows:
\begin{align*}
	\int_{Q_{2^{2 k}, m}} |\omega_k (x)| \, dx & = \Omega_k^+ |B(x_{k, m}, \delta_k)| -\Omega_k^{-} |\mathcal{C}_{k, m}| \\
	& = \pi \bigg[\Omega_k^+ \delta_k^2 + \frac{\Omega_k^+ \delta_k^2}{R_k^2-\delta_k} (R_k^2 - \delta_k) \bigg] \\
	& =  2\pi \Omega_k^+ \delta_k^2. 
\end{align*}
Then, under the choice of $\Omega_k^+$ given by
\begin{equation}\label{Omega_k^+}
	\Omega_k^+ = \frac{2^{-2 k}}{2 \pi \delta_k^2},
\end{equation}
we obtain that
\begin{equation}\label{EstimOC}
	\int_{Q_{2^{2 k}, m}} |\omega_k (x)| \, dx  = 2^{-2 k}
\end{equation}
and
$$
	\int_{Q_0}  |\omega_k (x)| \, dx = \sum_{m=1}^{2^{2 k}} \int_{Q_{2^{2 k}, m}} |\omega_k (x)| \, dx = 2^{2 k} 2^{-2 k} = 1. 
$$
Hence $\omega_k$ is $L^1$-normalized for every $k$. 

\subsection{$\{\omega_k\}$ is uniformly bounded in $M^{\log, 1}$.}\label{SectionMorreyLogComp}
We will make use of the well-known fact that Morrey norms can be characterized in terms of dyadic cubes, namely, $M^{\log, 1}$ is formed by all those  $\omega \in \mathcal{M}$ such that
\begin{equation}\label{EquivMorrey}
	 \sup_{\nu \in \N} \sup_{Q \in \mathbb{D}_\nu} |\log |Q|| |\omega|(Q) < \infty, 
\end{equation}
where $\mathbb{D}_\nu$ denotes the set of all dyadic cubes  with side length $2^{-\nu}$ contained in the unit cube. 

Next we shall use \eqref{EquivMorrey} to estimate $\|\omega_k\|_{M^{\log, 1}}$ for every $k \in \N$. Assume first that  $Q \in \mathbb{D}_\nu$ with $\nu > 2^{2 k}$. Then there exists a unique ancestor $\widehat{Q} \in \mathbb{D}_{2^{2 k}}$ such that $Q \subset \widehat{Q}$. If $\widehat{Q} \not \in \{Q_{2^{2 k}, m} : m= 1, \ldots, 2^{2 k}\}$, then we have $\omega_k (Q) =0$ since $\text{supp } \omega_k \subset E_k$ and $\widehat{Q} \cap E_k = \emptyset$. 

On the other hand, if $\widehat{Q} = Q_{2^{2 k}, m}$ for some $m \in \{1, \ldots, 2^{2 k}\}$, we are going to show that
\begin{equation}\label{6.0}
	\int_Q |\omega_k(x)| \, dx \lesssim |\log |Q||^{-1},
\end{equation}
where the hidden equivalence constants are independent of $Q$ and $k$. Indeed, we have (recalling \eqref{Omega_k-} and \eqref{Omega_k^+})
\begin{align}
	\int_Q |\omega_k(x)| \, dx  &= \int_{Q \cap (A_k \cup B_k) } |\omega_k(x)| \, dx \nonumber \\
	&= \int_{Q \cap A_k} |\omega_k(x)| \, dx + \int_{Q \cap B_k} |\omega_k(x)| \, dx \nonumber \\
	& =  \Omega_k^+ |Q \cap A_k| - \Omega_k^{-} |Q \cap B_k| \nonumber \\
	& =  \Omega_k^+ |Q \cap A_k| + \frac{\Omega_k^+ \delta_k^2}{R_k^2-\delta_k}  |Q \cap B_k|. \label{6.3}
\end{align}
Now we fix the value of $R_k$ to be
\begin{equation}\label{Rk}
	R_k = \sqrt{c \delta_k},
\end{equation}
where $c > 1$ is a fixed constant (note that, in particular, $R_k > \sqrt{\delta_k}$). Under this choice, \eqref{6.3} reads as
\begin{equation}\label{6.3n}
	\int_Q |\omega_k(x)| \, dx  =  \Omega_k^+ |Q \cap A_k| + \frac{\Omega_k^+ \delta_k}{c-1}  |Q \cap B_k|.
\end{equation}

In order to estimate \eqref{6.3n}, we need to specify the value of $\delta_k$, namely, we set
\begin{equation}\label{delta_k}
	\delta_k = \frac{l_k^2}{8 c}  = \frac{2^{-2^{2 k+1}}}{8 c}.
\end{equation}
We first observe that this choice is compatible with the standing assumption \eqref{Contain}  since (cf. \eqref{Rk}) 
$$2 R_k = 2 \sqrt{c \delta_k} =  \frac{l_k}{\sqrt{2}} < l_k.$$

We are now ready to deal with \eqref{6.3n}. We distinguish two possible cases. Firstly, if $|Q| \leq \delta_k^2$ (or in other words $2^{-\nu} \leq \delta_k$) then, by \eqref{6.3n}, \eqref{Omega_k^+} and taking into account that $\delta_k \to 0$ as $k \to \infty$ (cf. \eqref{delta_k}), 
\begin{align*}
		\int_Q |\omega_k(x)| \, dx  &\leq \Omega_k^+ |Q| \bigg(1 + \frac{\delta_k}{c-1} \bigg)   \lesssim  2^{-2 k} \delta_k^{-2} |Q| \\
		& = 2^{-2 k} \delta_k^{-2} |Q| (-\log |Q|) (-\log |Q|)^{-1} \\
		& \lesssim 2^{-2 k} \delta_k^{-2} \delta_k^2 (-\log \delta_k)  (-\log |Q|)^{-1} \\
		& \approx 2^{-2 k} (-\log l_k) (-\log |Q|)^{-1} \\
		& \approx 2^{-2 k} 2^{2 k} (-\log |Q|)^{-1}  = (-\log |Q|)^{-1}.
\end{align*}
Hence \eqref{6.0} holds provided that $|Q| \leq \delta_k^2$. 

Secondly,  suppose that $|Q| > \delta_k^2$. Recall that $Q \subset Q_{2^{2 k}, m}$ and then, in light of \eqref{EstimOC}, 
\begin{equation}\label{2nM}
	\int_Q |\omega_k(x)| \, dx \leq \int_{Q_{2^{2 k}, m}} |\omega_k(x)| \, dx = 2^{-2 k}. 
\end{equation}
On the other hand, monotonicity properties of $\log$ and the definition of $\delta_k$ given in \eqref{delta_k} allow us to estimate
$$
	(-\log |Q|)^{-1} \gtrsim  (-\log \delta_k)^{-1} = \bigg(-\log \frac{2^{-2^{2 k+1}}}{8 c} \bigg)^{-1} \approx 2^{-2 k}.
$$
Inserting this into \eqref{2nM}, we arrive at
$$
	\int_Q |\omega_k(x)| \, dx  \lesssim (-\log |Q|)^{-1},
$$
that is, \eqref{6.0} also holds if $|Q| > \delta_k^2$. 

So far, we have shown that, for every $k \in \N$, 
\begin{equation}\label{Summary}
	|\omega_k|(Q) \lesssim (-\log |Q|)^{-1} \qquad \text{if} \qquad Q \in \mathbb{D}_\nu, \quad \nu > 2^{2 k}. 
\end{equation}
Moreover, the case $\nu = 2^{2^{2 k}}$ in \eqref{Summary} is also covered by the assertion \eqref{EstimOC}. Then it remains to check the fulfilment of \eqref{Summary} under the assumption $Q \in \mathbb{D}_\nu$ with $\nu < 2^{2 k}$. Indeed, there exists a unique $l \in \{1, \ldots, k\}$ such that $2^{2(l-1)} < \nu \leq 2^{2 l}$. A simple geometrical argument shows that $Q$ contains at most $2^{ 2 k}/ 2^{2 l} = 2^{2(k-l)}$ cubes within the family $E_k$ (see Figure \ref{Figure2} for the case $l=k$ and iterate the process for $l < k$). As a consequence, it follows from \eqref{EstimOC} that
\begin{align*}
	\int_Q |\omega_k(x)| \, dx  & =  \sum_{m \in \{1, \ldots, 2^{2 k}\}: Q_{2^{2 k}, m} \subset Q} \int_{Q_{2^{2 k}, m}} |\omega_k(x)| \, dx \\
	& \leq 2^{-2 k} 2^{2(k-l)} = 2^{-2 l} \\
	& \leq \nu^{-1} \approx  (-\log |Q|)^{-1}.
\end{align*}
This extends \eqref{Summary} to the remaining case $Q \in \mathbb{D}_\nu, \, \nu < 2^{2 k}$ and completes the proof of the desired assertion
$$
	\|\omega_k\|_{M^{\log, 1}} \lesssim 1 \qquad \text{for every} \qquad k \in \N. 
$$

\subsection{Concentration-cancellation for $\{\omega_k\}$ via sparseness} 
We are going to show an extremely bad  behaviour of sparse numbers of $\{\omega_k\}$ in the sense that
\begin{equation}\label{Somegak}
	 s_{2^{2 N}}(\{\omega_k\}) \to \infty \qquad \text{as} \qquad N \to \infty, 
\end{equation}
suggesting a scenario where concentration-cancellation phenomenon arises (see Theorem \ref{ThmDM2}.)  

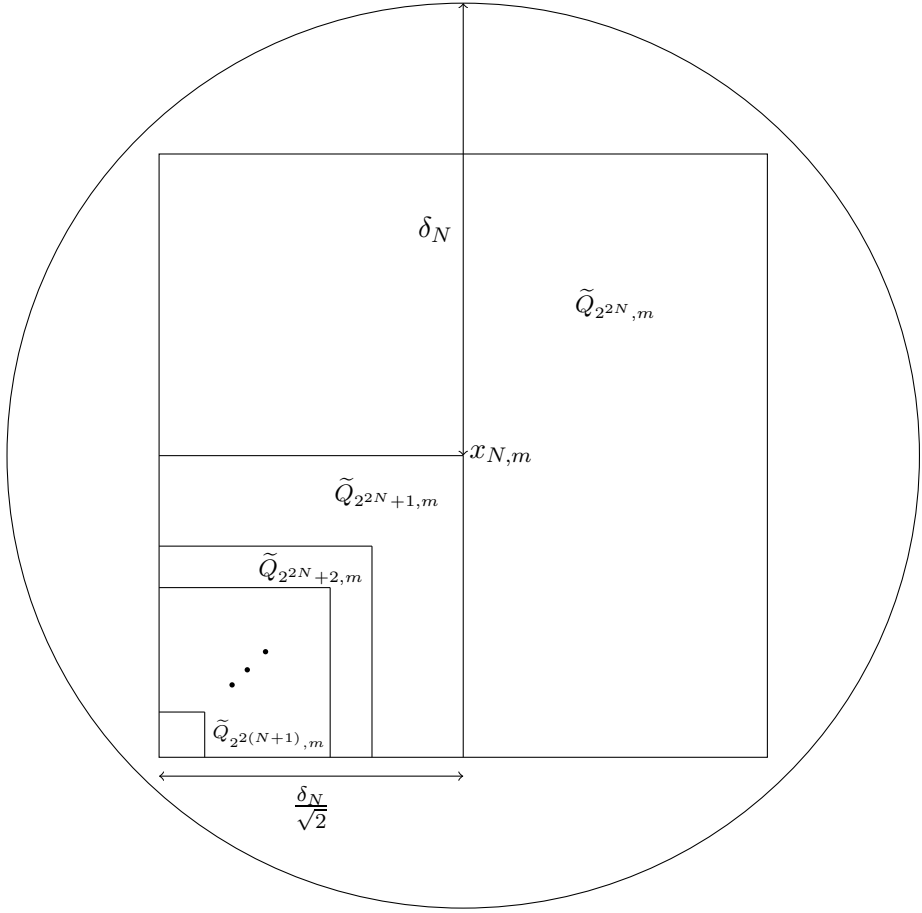
\begin{figure}[H]
\begin{center}

\begin{tikzpicture}

    % Dibujar el rectángulo grande
    \draw (0,0) rectangle (8,8);
    
    % Dibujar las líneas internas
    \draw (4,0) -- (4,4);
    \draw (0,4) -- (4,4);
    
        \draw (2.8,0) -- (2.8,2.8);
    \draw (0,2.8) -- (2.8,2.8);
    
        \draw (2.25,0) -- (2.25, 2.25);
    \draw (0,2.25) -- (2.25,2.25);
    
       \draw (0.6,0) -- (0.6,0.6);
    \draw (0,0.6) -- (0.6,0.6);
    
    \fill (0.96, 0.96) circle (1pt);
    
    \fill (1.4, 1.4) circle (1pt);

      \fill (1.16, 1.16) circle (1pt);
      
      \draw (4, 4) circle (6);

    % Etiquetas
    \node at (6,6) {\footnotesize $\widetilde{Q}_{2^{2 N},m}$};
        \node at (3,3.5) {\footnotesize $\widetilde{Q}_{2^{2N} + 1,m}$};
                \node at (2,2.5) {\footnotesize $\widetilde{Q}_{2^{2N} + 2,m}$};
                                \node at (1.45,0.3) {\tiny $\widetilde{Q}_{2^{2 (N+1)},m}$};

    \node at (4.5,4) { $x_{N,m}$};
    
    \draw[<->] (0, -0.25) -- (4, -0.25) node[midway, below] {$\frac{\delta_N}{\sqrt{2}}$};
        \draw[<->] (4, 4) -- (4, 10) node[midway, left] {$\delta_N$};

\end{tikzpicture}

\end{center}
\caption{\label{Figure4} Geometrical construction of sparse decompositions related to a vortex patch of $\omega_N$.}
\end{figure}

To show \eqref{Somegak} we propose a simple modification of the sparse decomposition constructed in Section \ref{Section4} adapted now to the geometry of $\omega_k$ (see Figure \ref{Figure3}). Recall that (see \eqref{SI1})
\begin{equation}\label{Cs0}
	s_{2^{2 N}} (\omega_N) = \sup_{S}  \bigg(\sum_{i \in I: \ell(Q_i) < 2^{-2^{2 N}}} |\omega_N|(Q_i)^2  \bigg)^{\frac{1}{2}},  \qquad N \in \N, 
\end{equation}
where the supremum runs over all sparse families of cubes. For any vortex patch of the vorticity $\omega_N$, say $Q_{2^{2 N}, m}$ for a fixed $m \in \{1, \ldots, 2^{2 N}\}$, we consider a cube $\widetilde{Q}_{2^{2 N}, m}$ such that
\begin{equation}\label{Cs1}
	\widetilde{Q}_{2^{2 N}, m} \subset B(x_{N, m}, \delta_N) \qquad \text{with} \qquad \ell(\widetilde{Q}_{2^{2 N}, m} ) \approx \delta_N. 
\end{equation}
To fix ideas, one may think that $\widetilde{Q}_{2^{2 N}, m}$ is the cube centered at $x_{N, m}$ with length of side $\delta_N/\sqrt{2}$. In particular, by \eqref{Cs1}, 
\begin{equation}\label{Cs3}
	\omega_N \equiv \Omega_N^+ \qquad \text{on} \qquad \widetilde{Q}_{2^{2 N}, m}. 
\end{equation}
Let us consider the dyadic\footnote{To avoid unnecessary technicalities, we may assume without loss of generality that all  cubes involved in this construction are dyadic.} tree of subcubes of  $\widetilde{Q}_{2^{2 N}, m}$  given by 
$$
	\mathscr{S}_m = \{\widetilde{Q}_{j, m} : j= 2^{2 N} +1, \ldots, 2^{2(N+1)} \},
$$
which is depicted in Figure \ref{Figure4}, and set
$$
	\mathscr{S} = \bigcup_{m=1}^{2^{2 N}} \mathscr{S}_m. 
$$
 Note that $\mathscr{S}_m$ (and hence $\mathscr{S}$ since $\mathscr{S}_m \cap \mathscr{S}_\ell = \emptyset$ if $m \neq \ell$) is a sparse family cubes that are not pairwise disjoint. This assertion was already justified in Subsection \ref{Section4n}. Furthermore, $\mathscr{S}$ is formed by cubes with sidelength at most $\frac{\delta_N}{\sqrt{2}} \approx 2^{-2^{2N+1}}$ (cf. \eqref{delta_k}).

 Applying \eqref{Cs0} and \eqref{Cs3} and taking into account the choices of $\Omega_N^+$ and $\delta_N$ given in  \eqref{Omega_k^+} and \eqref{delta_k}, respectively,  we can estimate
 \begin{align*}
 	s_{2^{2 N}} (\omega_N)  &\geq  \bigg(\sum_{Q \in \mathscr{S}} |\omega_N|(Q)^2  \bigg)^{\frac{1}{2}}  = \Omega_N^+ \bigg(\sum_{Q \in \mathscr{S}} |Q|^2  \bigg)^{\frac{1}{2}} \\
	& = \Omega_N^+ \bigg(\sum_{m=1}^{2^{2N}} \sum_{j= 2^{2 N} +1}^{2^{2(N+1)}} |\widetilde{Q}_{j, m}|^2 \bigg)^{\frac{1}{2}} \\
	%& = \Omega_N^+ \bigg(\sum_{m=1}^{2^{2N}} \sum_{j= 2^{2 N} +1}^{2^{2(N+1)}} 2^{-j 4} \bigg)^{\frac{1}{2}} \\
	& = \Omega_N^+ 2^N  \bigg( \sum_{j= 2^{2 N} +1}^{2^{2(N+1)}} 2^{-4j} \bigg)^{\frac{1}{2}}  \\
	& \approx \Omega_N^+ 2^N 2^{-2^{2 N} 2} \approx  2^{- N}  (\delta_N^{-1} 2^{-2^{2N}})^2 \\
	& \approx 2^{-N} (2^{2^{2N +1}} 2^{-2^{2N}})^2  = 2^{-N} 2^{2^{2N+1}}.
 \end{align*}
 In particular, we derive (cf. \eqref{SIn4})
 $$
 s_{2^{2 N}} (\{\omega_k\}) =	\sup_{k \in \N} s_{2^{2 N}} (\omega_k) \geq s_{2^{2 N}} (\omega_N) \gtrsim 2^{-N} (2^{2^{2N +1}} 2^{-2^{2N}})^2  = 2^{-N} 2^{2^{2N+1}}
 $$
 and then \eqref{Somegak} follows.

\subsection{Construction of velocity fields $\{u_k\}$ related to Cantor sets} We start by recalling the following well-known fact \cite[p. 308]{DiPernaMajda}: 
Consider any radial function $\omega(r), \, r=  |x| = \sqrt{x_1^2 + x_2^2}$, with $\text{supp } \omega \subset [0, R]$. Then $\omega_0$ defines a rotating eddy $u$, which is an exact solution of the steady  Euler equations \eqref{Euler} through the explicit formula
\begin{equation}\label{Repu}
	u (x) = r^{-2} x^\perp \int_0^r s \omega(s) \, ds,
\end{equation}
where $x^\perp = (-x_2, x_1)$. In particular, if we assume further that
\begin{equation}\label{AssSuppport}
	\int_0^\infty r \omega(r) \, dr = \int_0^R r \omega(r) \, dr = 0,
\end{equation}
then $\text{supp } u \subset \{x: |x| < R\}$. 

Given $k \in \N$, consider 
\begin{equation}\label{Defomegak}
	\omega^{k} (r) =  \left\{\begin{array}{cl}  \Omega_k^+ &\quad  \text{if} \quad  r \in [0, \delta_k), \\
	& \\
		\Omega_k^-  &\quad  \text{if} \quad r \in [\sqrt{\delta_k}, R_k), \\
		& \\
		0& \quad \text{otherwise}, 
		       \end{array}
                        \right.
\end{equation}
and the corresponding $u^{k}$ given by \eqref{Repu}. Elementary computations show that (cf. \eqref{Omega_k-})
$$
	\int_0^\infty r \omega^k(r) \, dr = \frac{\Omega_k^+ \delta_k^2 + \Omega_k^- (R_k^2 -\delta_k)}{2} =0;
$$
see also \eqref{ZeroVort}. Then \eqref{AssSuppport} holds and, in particular, the support of $u^{k}$ is contained in the ball  $\{x : |x| < R_k\}$. Note that $\omega_k$ defined in \eqref{DefVort} can be expressed in terms of adequate scalings of $\omega^k$ related to $x_{k, m}$, the set of centers of the cubes $Q_{2^{2 k}, m}$ involved in the definition of the Cantor set $E$ (cf. \eqref{Cantor}), namely, 
$$
	\omega_k (x) = \sum_{m=1}^{2^{2k}} \omega^k (|x-x_{k, m}|). 
$$
The corresponding sequence of velocity fields $u_k$ related to $E$ is given by
\begin{equation}\label{Defuk}
	u_k(x) = \sum_{m=1}^{2^{2 k}} u^k (x-x_{k, m}). 
\end{equation}
Note that the support of each term (i.e., $u^k (x-x_{k, m})$) in the definition of $u_k$  is $B(x_{k, m}, R_k)$ and, by construction, these supports are pairwise disjoint. Indeed, recall that the choice of $R_k$ given by \eqref{Rk} guarantees this fact, see also \eqref{Contain}. Then $u_k$ is the sum of steady solutions to  \eqref{Euler} that do not overlap and hence $u_k$ is also a steady solution to \eqref{Euler}. Moreover  (recall \eqref{Es})
\begin{equation}\label{Supportu}
	\text{supp } u_k \subset E_k. 
\end{equation}

\subsection{$\{u_k\}$ is uniformly bounded in $L^2$}\label{SectionBoundedL2}

Since the velocities $u^k (\cdot-x_{k, m})$ involved in the definition of $u_k$ (cf. \eqref{Defuk}) have pairwise disjoint supports, we have
\begin{equation*}
	\|u_k\|_{L^2}^2   = \sum_{m=1}^{2^{2 k}} \int_{B(x_{k, m}, R_k)} |u^k(x-x_{k, m})|^2 \, dx
\end{equation*}
and a simple change of variables leads to
\begin{equation}\label{L2normuk}
	\|u_k\|_{L^2}^2 = 2^{2 k} \int_{B(0, R_k)} |u^k(x)|^2 \, dx.
\end{equation}
Furthermore, thanks to \eqref{Repu} and \eqref{Defomegak}, $|u^k|$ can be explicitly computed in terms of the following radial function 
\begin{equation}\label{L2normuk2}
	|u^k (x)| = \frac{1}{r} \int_0^r s \omega^k (s) \, ds = \left\{\begin{array}{cl} \frac{r \Omega_k^+}{2} & \quad \text{if} \quad r \in [0, \delta_k), \\
	\frac{\delta_k^2 \Omega_k^+}{2 r} & \quad \text{if} \quad r \in [\delta_k, \sqrt{\delta_k}), \\
	\frac{\delta_k^2 \Omega_k^+ + \Omega_k^- (r^2 -\delta_k)}{2 r} & \quad \text{if} \quad r \in [\sqrt{\delta_k}, R_k), \\
	0 & \quad \text{if} \quad r \in [R_k, \infty).
	\end{array}
	\right.
\end{equation}
Applying polar coordinates to \eqref{L2normuk} and using \eqref{L2normuk2}, we obtain
\begin{align*}
	2^{-2 k}\|u_k\|_{L^2}^2  &= \int_0^{R_k} r |u^k(r)|^2 \, dr \\
	& = \bigg(\frac{\Omega_k^+}{2} \bigg)^2 \int_0^{\delta_k} r^3 \, dr + \bigg(\frac{\delta_k^2 \Omega_k^+}{2} \bigg)^2 \int_{\delta_k}^{\sqrt{\delta_k}} \frac{dr}{r}  \\
	& \hspace{2cm} + \frac{1}{4} \int_{\sqrt{\delta_k}}^{R_k}  \Big(\delta_k^2 \Omega_k^+ + \Omega_k^- (r^2 -\delta_k) \Big)^2 \, \frac{dr}{r} \\
	& = \frac{1}{4} \bigg(\frac{ \delta_k^2 \Omega_k^+}{2} \bigg)^2  +  \bigg(\frac{\delta_k^2 \Omega_k^+}{2} \bigg)^2 \frac{|\log \delta_k|}{2} + \bigg(\frac{\delta_k^2 \Omega_k^+}{2 (R_k^2-\delta_k)} \bigg)^2   \int_{\sqrt{\delta_k}}^{R_k}  (R_k^2 - r^2)^2 \, \frac{dr}{r} \\
	& \lesssim (\delta_k^2 \Omega_k^+)^2 + (\delta_k^2 \Omega_k^+)^2  |\log \delta_k| + \bigg(\frac{\delta_k^2 \Omega_k^+}{R_k^2-\delta_k} \bigg)^2 R_k^4 \log \bigg(\frac{R_k}{\sqrt{\delta_k}} \bigg) \\
	& \approx (\delta_k^2 \Omega_k^+)^2  \bigg[1 + |\log \delta_k| + \frac{R_k^4}{\delta_k^2} \bigg]  \\
	& \approx (\delta_k^2 \Omega_k^+)^2  |\log \delta_k| \approx 2^{-4 k} 2^{2 k} = 2^{-2k}, 
\end{align*}
where we have also used  \eqref{Omega_k-}, \eqref{Rk}, \eqref{Omega_k^+} and \eqref{delta_k}. The above computations  show that
$$
	\|u_k\|_{L^2} \lesssim 1 \qquad \forall k \in \N. 
$$

\subsection{$\{u_k\}$ converges weakly in $L^2$ to the trivial solution, but not strongly} Recall the well-known fact  that every bounded sequence in a Hilbert space admits a weakly convergent subsequence. In particular, by Section \ref{SectionBoundedL2}, the sequence of velocities $\{u_k\}$ given by \eqref{Defuk} is weakly convergent (possibly passing to a subsequence). Furthermore, we have (cf. \eqref{Contain} and \eqref{Es})
\begin{equation}\label{Pro1}
	\text{supp } u_k \subset \bigcup_{m=1}^{2^{2 k}} B(x_{k, m}, R_k) \subset \bigcup_{m=1}^{2^{2 k}} Q_{2^{2 k}, m} = E_k
\end{equation}
 with
 \begin{equation}\label{Pro2}
 |E_k| = \sum_{m=1}^{2^{2 k}} |Q_{2^{2 k}, m}| = \sum_{m=1}^{2^{2 k}} 2^{-2^{2k}} = 2^{-2^{2k}}  2^{2 k} \to 0 \qquad \text{as} \qquad k \to \infty.
 \end{equation}
 As a consequence of \eqref{Pro1} and using also H\"older's inequality, we derive for every $\varphi \in C^\infty_c$, 
 \begin{align*}
 	\bigg|\int_{\R^2} u_k(x) \varphi(x) \, dx \bigg| &\leq \int_{E_k} |u_k(x)| |\varphi(x)| \, dx \\
	&\leq \|u_k\|_{L^2} \bigg( \int_{E_k} |\varphi(x)|^2 \, dx \bigg)^{1/2} \\
	& \lesssim \|\varphi\|_{L^\infty} |E_k|^{1/2}.
\end{align*}
In view of \eqref{Pro2} and the basic fact that smooth functions with compact support are dense in $L^2$, we can conclude that 
$$
	u_k \rightharpoonup 0 \qquad \text{in} \qquad L^2. 
$$
 
On the other hand, $u_k \not \to 0$ in $L^2$. Indeed, computations carried out in Section \ref{SectionBoundedL2} show in particular that
\begin{align*}
	\|u_k\|_{L^2}^2  \geq  2^{2 k} \bigg(\frac{\delta_k^2 \Omega_k^+}{2} \bigg)^2 \int_{\delta_k}^{\sqrt{\delta_k}} \frac{dr}{r} \gtrsim 1. 
\end{align*}

\subsection{Characterization of the reduced defect measure $\theta$ in terms of the measure $\omega$}\label{RDF}
Let $Q$ be any dyadic cube in $\mathbb{D}_\nu$. Without loss of generality, we may assume that $k$ is large enough such that $2^{2 k} > \nu$.  We shall distinguish two possible cases: Firstly, if $Q \cap E_k = \emptyset$  then, by \eqref{Supportu},  $u_k$ vanishes on $Q$ and, in addition, $\delta_{x_{k, m}}(Q) =0$ for every $m=1, \ldots, 2^{2 k}$. Hence
\begin{equation}\label{Meas1}
	\int_Q |u_k(x)|^2 \, dx = 0 = \sum_{m=1}^{2^{2 k}} 2^{-2 k} \delta_{x_{k, m}}(Q).
\end{equation}

Secondly, suppose that $Q \cap E_k \neq \emptyset$. In particular, using the dyadic structure and the fact that $2^{2 k} > \nu$,  the sets
$$
	\mathcal{V}_k = \{Q_{2^{2 k}, m} : Q_{2^{2 k}, m} \subset Q\}
$$
are non-empty. Furthermore, it is plain to see that, for every fixed $k$ and $Q$,  the map
$
	Q_{2^{2 k}, m} \mapsto x_{k, m}
$
induces a bijection between $\mathcal{V}_k$ and 
$$
	\Lambda_k=	 \{x_{k, m} : x_{k, m} \in Q\}. 
$$

The energy of $u_k$ in $Q$ can be easily computed as
\begin{align*}
	\int_Q |u_k(x)|^2 \, dx &=  \sum_{m : Q_{2^{2 k}, m} \subset Q} \int_{B(x_{k, m}, R_k)} |u^k(x-x_{k, m})|^2 \, dx \\
	& = |\mathcal{V}_k| \int_{B(0, R_k)} |u^k(x)|^2 \, dx. 
\end{align*}
Furthermore, in Section \ref{SectionBoundedL2} we proved that $ \int_{B(0, R_k)} |u^k(x)|^2 \, dx \approx 2^{-2 k}$, which gives
$$
	\int_Q |u_k(x)|^2 \, dx \approx |\mathcal{V}_k| 2^{-2 k}.
$$
Based now on the bijection between the sets $\mathcal{V}_k$ and $\Lambda_k$, we can rewrite the latter as follows
\begin{equation}\label{Meas2}
	\int_Q |u_k(x)|^2 \, dx \approx |\Lambda_k| 2^{-2 k}  = \sum_{m : x_{k, m} \in Q} 2^{-2 k} = \sum_{m=1}^{2^{2 k}} 2^{-2 k} \delta_{x_{k, m}}(Q). 
\end{equation}

Recall the definition of $\omega$ given in Section \ref{Step2}. Putting together \eqref{Meas1} and \eqref{Meas2},  we obtain
$$
	\theta(Q) = \limsup_{k \to \infty} \int_Q |u_k(x)|^2 \, dx \approx \omega (Q)
$$
for all dyadic cubes $Q$. In particular,  this assertion can be extended to all open sets and characterizes those sets for which  strong convergence of $\{u_k\}$ to the weak limit (i.e., $0$) hold.

\subsection{Characterization of concentration sets}
By virtue of Section \ref{RDF}, concentration sets are characterized in terms of $\omega$, more precisely, 
\begin{equation*}
A \quad \text{is a concentration set for} \quad  \{u_k\} \iff \omega(A) > 0. 
\end{equation*}

We claim that the Cantor set $E$  is a concentration set for $\{u_k\}$ of dimension zero. Indeed, recall that (cf. \eqref{Cantor} and \eqref{nested})
$$
	E = \bigcap_{k =0}^\infty E_k \qquad \text{with} \qquad \qquad E_{k+1} \subset E_k.
$$
Furthermore, by Section \ref{RDF}, we have  
\begin{equation}\label{CS1}
\theta(E_k^c) = \omega(E_k^c)= 0.
\end{equation}
 It remains to show that $E$ has dimension zero: Given any $\delta > 0$ and $\gamma > 0$, we can always find $m \in \N$ large enough such that the following conditions
$$
	\sum_{k=m}^\infty l_k^\gamma 2^{2 k} \approx l_m^{\gamma} 2^{2 m} < \delta \qquad \text{and} \qquad \sup_{k \geq m} l_k = l_m < \delta 
$$
are satisfied simultaneously. Then, for this choice of $m$, we have
$$
	E \subset E_m = \bigcup_{k=m}^\infty E_k
$$
with $\theta (E_m^c) =0$ (see \eqref{CS1}). This proves that $E$ has dimension zero. \qed

\section*{Acknowledgements}
O. Dom\'inguez is supported by the AEI grant RYC2022-037402-I.  D. Spector is supported by the National Science and Technology Council of Taiwan under research grant numbers 110-2115-M-003-020-MY3/113-2115-M-003 -017 -MY3 and the Taiwan Ministry of Education under the Yushan Fellow Program.

\end{document}